\documentclass[12pt,a4paper,reqno]{amsart}
\usepackage{graphics,epic}
\usepackage{amsmath,amssymb, amsthm}
\usepackage[all,2cell]{xy}

\newtheorem{theorem}{Theorem}[section]
\newtheorem*{theorem*}{Theorem}

\newtheorem{proposition}[theorem]{Proposition}

\newtheorem*{conjecture*}{Conjecture}

\newtheorem{remark}[theorem]{Remark}
\newtheorem{definition}[theorem]{Definition}

\newcommand{\cc}{{\mathcal C}}

\newcommand{\cu}{{\mathcal U}}
\newcommand{\cv}{{\mathcal V}}

\renewcommand{\hat}[1]{\widehat{#1}}

\newcommand{\ot}{\otimes}

\newcommand{\id}{{\rm id}}
\newcommand{\im}{{\rm im}}
\newcommand{\Hom}{{\rm Hom}\,}

\newcommand{\End}{{\rm End}\,}
\newcommand{\Res}{{\rm Res}\,}

\newcommand{\Z}{\mathbb{Z}}

\newcommand{\C}{\mathbb{C}}

\def\Res{{\rm Res}}
\def\wt{{\rm wt}}

\def\C{{\mathbb C}}

\def\Z{{\mathbb Z}}

\def\Y{{\mathcal Y}}

\def\1{{\bf 1}}
\def\tr{{\rm tr}}
\def \End{{\rm End}}
\def \Hom{{\rm Hom}}

\def \pf{\noindent {\bf Proof: \,}}
\def\theequation{5.\arabic{equation}}

\def \h{\mathfrak{h}}
\def \l{\lambda}
\def \w{\omega}
\def \g{\mathfrak{g}}

\begin{document}

\title[Extensions of $L_{sl_2}(k,0)$ and preunitary VOAs with central charges $c<1$]{The extensions of $L_{sl_2}(k,0)$ and preunitary  vertex operator algebra with central charges $c<1$}

\author{ Chongying Dong}
\address{Chongying Dong, Department of Mathematics, University of
California\\ Santa Cruz, CA 95064}
\email{dong@ucsc.edu}
\thanks{The first author was supported by NSF grant DMS-1404741, NSA grant H98230-14-1-0118 and China NSF grant 11371261}
\author{Xingjun Lin}
\address{Xingjun Lin, Institute of Mathematics, Academia Sinica, Taipei 10617, Taiwan}
\email{linxingjun@math.sinica.edu.tw}
\begin{abstract}
The extensions of the affine vertex operator algebras $L_{sl_2}(k,0)$ and
the preunitary vertex operator algebras with central charges $c<1$ are classified. In particular, the unitary
vertex operator algebras with central charges $c<1$ are classified.
\end{abstract}
\maketitle
\section{Introduction \label{intro}}
\def\theequation{1.\arabic{equation}}
\setcounter{equation}{0}
One of important problems in conformal field theories is the classification of rational conformal field theory. It is well known that the mathematic structures underlying conformal field theory can be described by conformal nets and by vertex operator algebras. The conformal nets with central charges $c<1$ were classified in \cite{KL}. Rational vertex operator algebras with central charge $c<1$ were partly classified in  \cite{DZ}, \cite{M} and rational  vertex operator algebras with central charge $c=1$ were partially classified in \cite{DM1}, \cite{DJ1}, \cite{DJ2}, \cite{DJ3}. In this paper, we classify the extensions of affine vertex operator algebra $L_{sl_2}(k,0)$ associated
to $sl_2$ for any $k\geq 0$ and the preunitary vertex operator algebra central charges $c<1.$ A preunitary vertex operator algebra is  a vertex operator algebra which is also a unitary module for the Virasoro algebra. Consequently, the unitary vertex operator algebras \cite{DLin} with central charges $c<1$ are classified.

The results in this paper are heavily influenced by the modular invariants of vertex operator algebra $L_{sl_2}(k,0)$
and the Virasoro vertex operator algebra $L(c_m,0)$ given in \cite{CIZ1}, \cite{CIZ2}, \cite{Ka}. The key fact for the  classification of  preunitary vertex operator algebras is that any  preunitary vertex operator algebra $V$ with central charge $c<1$ is an extension of the Virasoro vertex operator algebra $L(c_m, 0)$. It is well known that both
$L_{sl_2}(k,0)$ and $L(c_m,0)$ are rational and $C_2$-cofinite \cite{FZ}, \cite{Z}, \cite{DMZ}, \cite{W}.  By  results in \cite{HKL} and \cite{ABD}, we know that any extension $V$ of $L_{sl_2}(k,0)$ or $L(c_m,0)$  is rational and $C_2$-cofinite. Thus there is a modular invariant of $L_{sl_2}(k,0)$ or $L(c_m, 0)$ corresponding to $V$ \cite{DLN}. On the other hand, the modular invariants of $L_{sl_2}(k,0)$ and $L(c_m, 0)$ were classified in \cite{CIZ1}, \cite{CIZ2} and \cite{Ka}. We will construct vertex operator algebras realizing the modular invariants of $L_{sl_2}(k,0)$ and $L(c_m, 0)$ and further classify the extensions of $L_{sl_2}(k,0)$ and $L(c_m, 0).$

 The main problem is to construct vertex operator algebras realizing the modular invariants of $L_{sl_2}(k,0)$ and $L(c_m, 0)$.  Some  simple current extensions of $L_{sl_2}(k,0)$ and  $L(c_{m}, 0)$ have been constructed in \cite{Li} and \cite{LLY}. We will construct the other vertex operator algebras by using conformal embedding, coset construction and mirror extension. Although the classification of extensions of $L_{sl_2}(k,0)$ is independent
 of the classification of extensions of $L(c_{m},0),$ it plays essential role in the classification of extensions of $L(c_{m},0).$ This is why we do the classification of extensions of $L_{sl_2}(k,0)$ first in this paper.

 We should mention that not every modular invariant given in \cite{CIZ1}, \cite{CIZ2}, \cite{Ka} can be realized by vertex operator algebra extensions  of $L_{sl_2}(k,0)$ and  $L(c_{m}, 0).$ It seems that the other modular invariants can be realized by vertex operator superalgebra extensions. For example, the
 modular invariants of type $D_{2n+1}$ associated to $L_{sl_2}(4n-2,0)$ are realized by  the vertex operator superalgebra extension   of
 $L_{sl_2}(4n-2,0)$ \cite{DH}. We will deal with the other modular invariants using the vertex operator superalgebras in another paper.

The paper is organized as follows: In Section 2, we recall some basic definitions and facts of vertex operator algebras. In Section 3, we construct vertex operator algebras realizing the modular invariants of $L_{sl_2}(k,0)$ and further classify rational and $C_2$-cofinite extension vertex operator algebras of $L_{{sl_2}}(k,0)$. In Section 4, we construct vertex operator algebras realizing the modular invariants of  Virasoro vertex operator algebras $L(c_{m}, 0)$ and then classify preunitary vertex operator algebras with central charge $c<1$.

\section{Preliminaries }
\def\theequation{1.\arabic{equation}}
\setcounter{equation}{0}
\subsection{Vertex operator algebras}In this subsection we briefly review  some basic notions and facts in the theory of vertex operator algebras from \cite{DLM}, \cite{FHL}, \cite{FLM},  \cite{LL} and \cite{Z}. Let $(V, Y, \1, \w)$ be a vertex operator algebra as defined in \cite{FLM} (see also \cite{B}). $V$ is called of {\em CFT} type if $V=\bigoplus_{n\geq 0}V_n$ and $\dim V_0=1$.

 An anti-linear  map $\phi:V\to V$ is called an {\em anti-linear automorphism} of $V$ if $\phi(\bf1)=\bf1, \phi(\omega)=\omega$ and
$\phi(u_nv)=\phi(u)_n\phi(v)$ for any $ u, v\in V$ and $n\in
\mathbb{Z}$.
Let $(V, Y, \1, \w)$ be a vertex operator algebra and $\phi: V\to V$ be an {\em anti-linear involution}, i.e. an anti-linear automorphism of order $2$.  $(V, \phi)$ is called {\em unitary} \cite{DLin} if there exists a positive definite Hermitian form $(,): V\times V\to \mathbb{C}$
which is $\C$-linear on the first vector and anti-$\C$-linear on the second vector such that the following invariant property holds: For any $a, u, v\in V$\\
 $$(Y(e^{zL(1)}(-z^{-2})^{L(0)}a, z^{-1})u, v)=(u,
Y(\phi(a), z)v),$$where $L(n)$ is defined by $Y(\w, z)=\sum_{n\in \Z}L(n)z^{-n-2}$.

A vertex operator algebra $V$ is called {\em preunitary} if it is a unitary module for the Virasoro algebra spanned by
$\{L(n), c|n\in\Z\}.$ It is clear that a unitary vertex operator algebra is preunitary.

Let $V$ be a vertex operator algebra. A {\em weak  $V$-module} $M$ is a vector space equipped
with a linear map
\begin{align*}
Y_{M}:V&\to (\End M)[[x, x^{-1}]],\\
v&\mapsto Y_{M}(v,x)=\sum_{n\in\Z}v_nx^{-n-1},\,v_n\in \End M
\end{align*}
satisfying the following conditions: For any $u\in V,\ v\in V,\ w\in M$ and $n\in \Z$,
\begin{align*}
&\ \ \ \ \ \ \ \ \ \ \ \ \ \ \ \ \ \ \ \ \ \ \ \ \ \ u_nw=0 \text{ for } n>>0;\\
&\ \ \ \ \ \ \ \ \ \ \ \ \ \ \ \ \ \ \ \ \ \ \ \ \ \ Y_M(\1, x)=\id_M;\\
\begin{split}
&x_{0}^{-1}\delta\left(\frac{x_{1}-x_{2}}{x_{0}}\right)Y_{M}(u,x_{1})Y_M(v,x_{2})-x_{0}^{-1}\delta\left(
\frac{x_{2}-x_{1}}{-x_{0}}\right)Y_M(v,x_{2})Y_M(u,x_{1})\\
&\quad=x_{2}^{-1}\delta\left(\frac{x_{1}-x_{0}}{x_{2}}\right)Y_M(Y(u,x_{0})v,x_{2}).
\end{split}
\end{align*}

A weak
 $V$-module  $M$ is called an \textit{admissible $V$-module} if $M$ has a $\Z_{\geq
0}$-gradation $M=\bigoplus_{n\in\Z_{\geq 0}}M(n)$ such
that
\begin{align*}\label{AD1}
a_mM(n)\subset M(\wt{a}+n-m-1)
\end{align*}
for any homogeneous $a\in V$ and $m,\,n\in\Z$.
An admissible $V$-module $M$ is said to be
\textit{irreducible} if $M$ has no non-trivial admissible weak
$V$-submodule. When an admissible $V$-module $M$ is a
direct sum of irreducible admissible submodules, $M$ is called
\textit{completely reducible}.
A vertex operator algebra $V$ is said to be \textit{rational} if
any  admissible $V$-module is completely reducible. 
It was proved in \cite{DLM1} that if $V$ is rational then there are only finitely many irreducible admissible $V$-modules up to isomorphism.

A {\em  $V$-module} is a weak $V$-module $M$ which carries a $\C$-grading induced by the spectrum of $L(0)$, that is  $M=\bigoplus_{\lambda\in\C}
M_{\lambda}$ where
$M_\lambda=\{w\in M|L(0)w=\lambda w\}$. Moreover one requires that $M_\lambda$ is
finite dimensional and for fixed $\lambda\in\C$, $M_{\lambda+n}=0$
for sufficiently small integer $n$.

If $(V, \phi)$ is a unitary vertex operator algebra, a $V$-module $(M, Y_M)$  is called a {\em unitary}
 $V$-module if there exists a positive definite
Hermitian form $(,)_M: M\times M\to \mathbb{C}$ which is $\C$-linear on the first vector and anti-$\C$-linear on the second vector such that the following invariant property:
$$(Y_M(e^{zL(1)}(-z^{-2})^{L(0)}a, z^{-1})w_1, w_2)_M=(w_1, Y_M(\phi(a),
z)w_2)_M$$ holds for $a \in V $ and $w_1, w_2\in M$.

Let $M = \bigoplus_{\lambda\in \mathbb{C}}{M_{\lambda}}$ be a $V$-module. Set $M'
= \bigoplus_{\lambda \in \mathbb{C}}{M_\lambda^*}$, the restricted
dual of $M$. It was proved in \cite{FHL} that $M'$ is naturally a
$V$-module where the vertex operator map denoted by $Y'$ is defined
by the property
$$\langle Y'(a, z)u', v\rangle  = \langle u', Y(e^{zL(1)}(-z^{-2})^{L(0)}a, z^{-1})v\rangle ,$$\\for $a\in V, u'\in
M'$ and $v\in M$. The $V$-module $M'$ is called the {\em contragredient
module} of $M$. It was also proved in \cite{FHL} that if $M$ is irreducible,  then so
is $M'$, and that $(M')'\simeq M$. A $V$-module $M$ is called {\em self dual} if $M\cong M'$.
\subsection{Modular invariance of trace functions}
We now turn our discussion to the modular invariance property in the theory of vertex operator algebras.  Let $V$ be a rational  vertex operator algebra, $M^0=V, M^1,...,M^p$ be all the irreducible $V$-modules. Then $M^i, 0\leq i \leq p$, has the form $$M^i=\bigoplus_{n=0}^{\infty}M^i_{\lambda_i+n},$$ with $M^i_{\lambda_i}\neq 0$ for some number $\lambda_i$ which is called {\em conformal weight} of $M^i$.
 Let $\h =\{\tau\in \mathbb{C}| \im\tau>0\}$, for any irreducible $V$-module $M^i$ the trace function associated to $M^i$ is defined as follows: For any homogenous element $v\in V$ and $\tau\in \h$,
\begin{equation*}
Z_{M^i}(v,\tau):=\tr_{M^i}o(v)q^{L(0)-c/24}=q^{\lambda_i-c/24}\sum_{n\in\mathbb{Z}^+} \tr_{M^i_{\l_i+n}}o(v)q^n,
\end{equation*}
where $o(v)=v(\wt v-1)$ and $q=e^{2\pi i\tau}$. Recall  that a vertex operator algebra $V$ is called {\em $C_2$-cofinite} if $\dim V/C_2(V)<\infty$, where $C_2(V)=\langle u_{-2}v|u, v\in V \rangle$. Then $Z_{M^i}(v,\tau)$ converges to a holomorphic function on the domain $|q| < 1$ if $V$ is $C_2$-cofinite \cite{DLM2}, \cite{Z}.
Recall that the full modular group $SL(2, \mathbb{Z})$ has generators $S=\left(\begin{array}{cc}0 & -1\\ 1 & 0\end{array}\right)$, $T=\left(\begin{array}{cc}1 & 1\\ 0 & 1\end{array}\right)$ and acts on $\h$ as follows:$$\gamma: \tau\longmapsto \frac{a\tau+b}{c\tau+d}, \  \gamma=\left(\begin{tabular}{cc}
$a$ $b$\\
$c$ $d$\\
\end{tabular}\right) \in SL(2, \mathbb{Z}).$$
The following theorem was proved in \cite{Z} (also see \cite{DLM2}).
\begin{theorem}\label{minvariance}
 Let $V$ be a rational and $C_2$-cofinite vertex operator algebra with the irreducible $V$-modules $M^0,...,M^p.$   Then the vector space spanned by $Z_{M^0}(v,\tau),..., Z_{M^p}(v,\tau)$ is invariant under the action of $SL(2, \mathbb{Z})$ defined above, i.e. there is a representation $\rho$ of $SL(2, \mathbb{Z})$ on this vector space and the transformation matrices are independent of the choice of $v\in V$.
\end{theorem}

 We now assume that $V$ is a vertex operator algebra satisfying the following conditions:

(i) $V$ is a simple $CFT$ type vertex operator algebra  and is self dual;

(ii) $V$ is rational and $C_{2}$-cofinite.\\
Let $M^0, ..., M^p$ be all the $V$-irreducible modules.  Set $${\bf Z}(u,\tau)=(Z_{M^0}(u,\tau), ..., Z_{M^p}(u,\tau))^T,$$ then we have the following fact proved in \cite{DLN}.
\begin{proposition}\label{unique}
 If $A=(a_{ij})$ is a  matrix such that  for any $u,v\in V$,$${\bf Z}(u,\tau_1)^TA\overline{{\bf Z}(v,\tau_2)}=0,$$ then $A=0.$
\end{proposition}
\subsection{Modular invariants of vertex operator algebras}
In this subsection we recall some facts about the modular invariants of vertex operator algebras.
\begin{definition}
{\rm
Let $V$ be a vertex operator algebra satisfying  conditions (i) and (ii), $M^0=V, M^1,...,M^p$ be all the irreducible $V$-modules. A {\em modular invariant} of $V$ is a $(p+1)\times (p+1)$-matrix $X$ satisfying the following conditions:

(M1) The entries of $X$ are nonnegative integers;

(M2) $X_{00}=1;$

(M3) $XS=S X$ and $XT=TX$, where we use $S, T$ to denote the modular transformation matrix $\rho(S)$ and $\rho(T)$ respectively.
}
\end{definition}

  In the following we shall define a modular invariant of $V$ associated to  an extension vertex operator algebra of $V$. Recall that a vertex operator algebra $U$ is called an {\em extension vertex operator algebra} of $V$ if $V$ is a vertex operator subalgebra of $U$ and $V$, $U$ have the same conformal vector.
    \begin{theorem}\label{C 2}\cite{ABD}
Let $V$ be a $C_2$-cofinite vertex operator algebra and $U$ be an extension vertex operator algebra of $V$. Then $U$ is $C_2$-cofinite.
\end{theorem}
The following fact was essentially proved  in \cite {HKL}.
  \begin{theorem}\label{cvoa2}
Let $V$ be a vertex operator algebra satisfying  conditions $(i)$ and $(ii)$, $U$ be a simple extension vertex operator algebra of $V$. Then $U$ is rational.
\end{theorem}

  We now assume that $V$ is a vertex operator algebra satisfying conditions (i) and
(ii) and $U$ is an extension vertex operator algebra of $V$ satisfying conditions (i) and
(ii).
For $u,v\in V$,  set
$$f_V(u,v,\tau_1,\tau_2)=\sum_{i=0}^pZ_i(u,\tau_1)\overline{Z_i(v,\tau_2)},$$ where $\tau_1, \tau_2 \in\mathfrak{h}$.

 Similarly, for $u,v\in U$, set
$$f_U(u,v,\tau_1,\tau_2)=\sum_{M}Z_M(u,\tau_1)\overline{Z_M(v,\tau_2)},$$
 where $M$ ranges through the equivalent classes of irreducible $U$-modules.
Since each irreducible $U$-module $M$ is a direct sum of irreducible $V$-modules, there exists a matrix $X=(X_{i,j})$ such that  $X_{ij}\geq 0$ for all $i,j$ and for $u,v\in V$,
$$f_U(u,v,\tau_1,\tau_2)=\sum_{i,j=0}^pX_{ij}Z_i(u,\tau_1)\overline{Z_i(v,\tau_2)}.$$
By Proposition \ref{unique} the matrix $X=(X_{i,j})$ is uniquely determined. By \cite{DLN} we have
\begin{theorem} \label{invariant}
The matrix $X$ is a modular invariant of $V$.
\end{theorem}
\subsection{Tensor category of $V$-modules}
In this subsection we review some facts about tensor category of $V$-modules. We refer the reader to \cite{BK} for the basic notions of tensor categories.
First, we  recall the  notions of intertwining operators and fusion
rules from \cite{FHL}.
Let $M^1$, $M^2$, $M^3$ be weak $V$-modules. An {\em intertwining
operator} $\mathcal {Y}$ of type $\left(\begin{tabular}{c}
$M^3$\\
$M^1$ $M^2$\\
\end{tabular}\right)$ is a linear map
\begin{align*}
\mathcal
{Y}: M^1&\rightarrow \Hom(M^2, M^3)\{x\},\\
 w^1&\mapsto\mathcal {Y}(w^1, x) = \sum_{n\in \mathbb{C}}{w_n^1x^{-n-1}}
\end{align*}
satisfying the following conditions: For any $v\in V, w^1\in M^1, w^2\in M^2$ and $\lambda \in \mathbb{C},$
\begin{align*}
&  w_{n+\lambda}^1w^2 = 0 \text{ for }n>>0;\\
&  \dfrac{d}{dx}\mathcal{Y}(w^1,
x)=\mathcal{Y}(L(-1)w^1, x);\\
\begin{split}
x_0^{-1}\delta(\frac{x_1-x_2}{x_0})Y_{M^3}(v, x_1)&\mathcal
{Y}(w^1, x_2)-x_0^{-1}\delta(\frac{x_2-x_1}{-x_0})\mathcal{Y}(w^1,
x_2)Y_{M^2}(v, x_1)\\
&=x_2^{-1}\delta(\frac{x_1-x_0}{x_2})\mathcal{Y}(Y_{M^1}(v, x_0)w^1, x_2).
\end{split}
\end{align*}

Denote the vector space of intertwining operators of type $\left(\begin{tabular}{c}
$M^3$\\
$M^1$ $M^2$\\
\end{tabular}\right)$ by $\mathcal{V}_{M^1,M^2}^{M^3}$.  The dimension  of $\mathcal{V}_{M^1,M^2}^{M^3}$ is called the
{\em fusion rule} for $M^1$, $M^2$ and $M^3$, and is denoted by $N_{M^1,M^2}^{M^3}$.
Assume that $V$ is
a  vertex operator algebra satisfying conditions (i) and (ii).  By Lemma
4.1 in \cite{H4}, one knows that for 
$u^{i_k}\in M^{i_k}$, 
\begin{align*}
\langle u^4
,\Y_{M^{i_1},M^{i_5}}^{M^{i_4}}(u^1 , z_1)\Y_{M^{i_2},M^{i_3}}^{M^{i_5}}(u^2 , z_2)u^3
\rangle,
\end{align*}
\begin{align*}
\langle u^4
,\Y_{M^{i_2},M^{i_6}}^{M^{i_4}}(u^1 , z_1)\Y_{M^{i_1},M^{i_3}}^{M^{i_6}}(u^2 , z_2)u^3
\rangle,
\end{align*}
are analytic on $|z_1| > |z_2| > 0$ and $|z_2| > |z_1| > 0$ respectively and can both be
analytically extended to multi-valued analytic functions on
$$R = \{(z_1, z_2) \in \C^2|z_1, z_2 \neq 0, z_1 \neq z_2\}.$$
We can lift the multi-valued analytic functions on $R$ to single-valued analytic functions
on the universal covering $\tilde{R}$ of $R$ as in \cite{H5}. We use
\begin{align*}
E\langle u^4
,\Y_{M^{i_1},M^{i_5}}^{M^{i_4}}(u^1 , z_1)\Y_{M^{i_2},M^{i_3}}^{M^{i_5}}(u^2 , z_2)u^3
\rangle,
\end{align*}
and
\begin{align*}
E\langle u^4
,\Y_{M^{i_2},M^{i_6}}^{M^{i_4}}(u^1 , z_1)\Y_{M^{i_1},M^{i_3}}^{M^{i_6}}(u^2 , z_2)u^3
\rangle,
\end{align*}
to denote the analytic functions respectively.
Let $\{\Y_{M^{i_1},M^{i_2}}^{M^{i_3};j}
|1\leq j\leq N_{M^{i_1},M^{i_2}}^{M^{i_3}}\}$ be a basis of $\cv_{M^{i_1},M^{i_2}}^{M^{i_3}}$. The linearly independency of
\begin{align*}
\{E\langle u^4
,\Y_{M^{i_1},M^{i_5}}^{M^{i_4};j_1}(u^1 , z_1)\Y_{M^{i_2},M^{i_3}}^{M^{i_5};j_2}(u^2 , z_2)u^3
\rangle| 1\leq j_1\leq N_{M^{i_1},M^{i_5}}^{M^{i_4}},  1\leq j_2\leq N_{M^{i_2},M^{i_3}}^{M^{i_5}}\}
\end{align*}
follows from \cite{H5}. Moreover, for any $M^1, M^3, M^3, M^4$,
\begin{align*}
&span\{E\langle u^4
,\Y_{M^{i_1},M^{i_5}}^{M^{i_4};j_1}(u^1 , z_1)\Y_{M^{i_2},M^{i_3}}^{M^{i_5};j_2}(u^2 , z_2)u^3
\rangle|j_1,  j_2,M^{i_5}\}\\
&=span\{E\langle u^4
,\Y_{M^{i_2},M^{i_6}}^{M^{i_4};j_1}(u^1 , z_1)\Y_{M^{i_1},M^{i_3}}^{M^{i_6};j_2}(u^2 , z_2)u^3
\rangle|j_1,  j_2, M^{i_6}\}.
\end{align*}

\vskip0.5cm
We now turn our discussion to tensor categories of $V$-modules.
 Let $M^1$, $M^2$ be weak $V$-modules. A {\em tensor product} for the ordered pair $(M^1, M^2)$ is a pair $(M, \Y^{12})$ consisting of a weak $V$-module $M$ and an intertwining operator $\Y^{12}$ of type $\left(\begin{tabular}{c}
$M$\\
$M^1$ $M^2$\\
\end{tabular}\right)$ satisfying the following universal property: For any weak $V$-module $W$ and any intertwining operator $I$ of type $\left(\begin{tabular}{c}
$W$\\
$M^1$ $M^2$\\
\end{tabular}\right)$, there exists a unique $V$-homomorphism $\psi$ from $M$ to $W$ such that $I=\psi\circ \Y^{12}$.

Assume that $V$ is a vertex operator algebra satisfying conditions (i) and (ii). The tensor product $\boxtimes$ of $V$-modules was defined in a series of papers \cite{HL1}, \cite{HL2}, \cite{HL3}, \cite{HL4} and \cite{H1}. For any $V$-module $W$, set $\theta_{W}=e^{2\pi i L(0)}.$ Then the $V$-module category $\cc_V$ is a modular tensor category (see  \cite{H6}).
\begin{theorem}\label{kvoa1}
 Let $V$ be a vertex operator algebra satisfying conditions $(1)$ and $(2)$. Then the $V$-module category $(\cc_V, \boxtimes)$ is a modular tensor category such that the dual of $W$ is isomorphic to $W'$.
 \end{theorem}
To describe the next results, we now recall from \cite{KO} the notion of associative commutative algebra in modular tensor category.
 \begin{definition}
{\rm  Let $\cc$ be a modular tensor category, $\1_{\cc}$ be the unit object in $\cc$. An {\em associative commutative algebra} $A$ in $\mathcal{C}$ is an object $A\in \mathcal{C}$ along with morphisms $\mu: A\otimes A\to A$ and $\iota_A:\1_{\cc}\hookrightarrow A$ such that the following conditions hold:

1. (Associativity) Compositions $\mu\circ (\mu\otimes \id)\circ a,\ \mu\circ (\id\otimes \mu):A\otimes (A\otimes A)\to A$ are equal, where $a$ denotes the associativity isomorphism $a: A\otimes (A\otimes A)\to (A\otimes A)\otimes A$.

2. (Commutativity) Composition $\mu \circ c_{AA}:A\otimes A\to A$ is equal to $\mu$, where $c_{AA}$ denotes the commutativity isomorphism.

3. (Unit) Composition $\mu \circ (\iota_A\otimes A):A=\1_{\cc}\otimes A \to A$ is equal to $\id_A$.

4. (Uniqueness of unit) $\dim \Hom _{\mathcal{C}}(\1_{\cc}, A)=1$.
}
\end{definition}
The following results were established in \cite{HKL}.
\begin{theorem}\label{kvoa2} If $U$ is an extension vertex operator algebra of $V$, then
$U$ induces an associative commutative algebra $A_{U}$ in $\cc_V$. Conversely, if $U$ is a $V$-module having integral conformal weight and $U$ is an associative commutative algebra  in $\cc_V$, then $U$ has a vertex operator algebra structure such that $U$ is an extension vertex operator algebra of $V$.
\end{theorem}
\subsection{Mirror extension of vertex operator algebras}
In this subsection we recall the mirror extensions of vertex operator algebras from \cite{Mirror}. Let $(U, Y, \1, \w)$ be a vertex operator algebra and
$(V, Y, \1, \w')$ be a vertex operator subalgebra of $U$ such
that $\w'\in U_2$ and $L(1)\w'=0$. Consider the space $V^c=\{v\in U|L'(-1)v=0\}$, it was proved in \cite{FZ} that $(V^c, Y, \1, \w-\w')$ is a vertex operator subalgebra of
$U$, the so-called the commutant vertex operator algebra of $V$ in $U$.
Assume further that  $(V, Y, \1, \w)$, $(V^c, Y, \1, \w-\w')$,
$(U, Y, \1, \w')$ satisfy conditions (i), (ii) and $(V^c)^c=V$. Then it was proved in \cite{Mirror} that $U$ has the following decomposition as $V\ot V^c$-module:$$U=V\ot V^c\oplus (\oplus_{i=1}^nM^i\ot N^i),$$ where $\{M^i|1\leq i\leq n\}$ $($resp. $\{N^i|1\leq i\leq n\}$$)$ are inequivalent irreducible $V$-modules $($resp. $V^c$-modules$)$. Let $V^e=V\oplus (\oplus_{i=1}^nm_iM_i)$ and $ (V^c)^e=V^c\oplus (\oplus_{i=1}^nm_iN_i'),$  where $m_i$'s are nonnegative integers. Then we have the following result  proved in \cite{Mirror}.
\begin{theorem}\label{kvoa3}
Assume that there is a vertex operator algebra structure on $V^e$ such that $V^e$ is an extension vertex operator algebra of $V$. Then $(V^c)^e$ has a vertex operator algebra structure such that $(V^c)^e$ is an extension vertex operator algebra of $V^c$. Moreover, $(V^c)^e$ is a simple vertex operator algebra if $V^e$ is a simple vertex operator algebra.
\end{theorem}

\section{Extension vertex operator algebras of  $L_{{sl_2}}(k,0)$ }
\def\theequation{3.\arabic{equation}}
\setcounter{equation}{0}

 In this section, we classify  extension vertex operator algebras of affine vertex operator algebras $L_{{sl_2}}(k,0)$. First, we recall some facts about vertex operator algebras associated to  affine Lie algebras from \cite{FZ} ,\cite{DL} and \cite{K}.  Let
$\g$ be a finite dimensional simple Lie algebra and $\h$ be a Cartan subalgebra of $\g$. Fix
a non-degenerate symmetric invariant bilinear form $(,)$ on $\g$ so that $(\theta, \theta)=2$, where $\theta$ is the
maximal root of $\g$. Consider the affine Lie algebra
$\hat{\g}=\g\otimes\C[t, t^{-1}]\oplus \C K$ with the commutation relations
$$[x(m), y(n)]=[x, y](m+n)+Km(x, y)\delta_{m+n, 0},$$
$$[\hat{\g}, K]=0,$$ where $x(m)=x\otimes t^{n}$.

For a complex number $k\in \C$, set
$$V_\g(k)=U(\hat{\g})/J_k,$$ where $J_k$ is the left ideal of
$U(\hat{\g})$ generated by $x(n)$ and $K-k$ for $x\in \g$ and $n\geq 0.$ It is
well-known that  $V_\g(k)$ has a vertex operator algebra structure if
$k\neq -h^{\vee}$, where $h^{\vee}$ is the dual Coxeter number of
$\g$ (see \cite{FZ}). Moreover, $V_\g(k)$ has a unique maximal proper $\hat{\g}$-submodule
$J(k)$ and  $L_\g(k, 0)=V_\g(k)/J(k)$ is a simple vertex operator algebra. As usual
we denote the corresponding irreducible highest weight module for $\hat{\g}$ associated to a highest weight $\lambda\in\h^*$
of $\g$ by $L_\g(k, \lambda)$. It was proved in \cite{FZ}, \cite{DLM2} that $L_\g(k, 0)$ 
satisfies conditions (i) and (ii) if
$k\in \Z^{+}$ and
\begin{align*}
\{L_\g(k, \lambda)|(\lambda, \theta)\leq k, \lambda\in \h^*\text{ is integral dominant}\}
 \end{align*}are the complete list of inequivalent irreducible $L_\g(k, 0)$-modules.

In the following we take $g=sl_2$. From the discussion above, we have $\{L_{sl_2}(k, j)| 0\leq j\leq k\}$ are the complete list of inequivalent irreducible $L_{{sl_2}}(k,0)$-modules. The modular invariants of vertex operator algebras $L_{{sl_2}}(k,0)$ were classified in \cite{CIZ1}, \cite{CIZ2} and \cite{Ka}:
\begin{theorem}\label{CIZ1}
 Let $k$ be a nonnegative integer. We use $Z_j(u,\tau)$ to denote the trace function associated to $L_{sl_2}(k,0)$-module $L_{sl_2}(k,j)$. Then any modular invariant of $L_{{sl_2}}(k,0)$ is equal to one of the following modular invariants:
{\tiny
$$
\begin{tabular}{|c c c|}
\hline

any k  & $\sum\limits_{1\leq j \leq k}Z_{j}(u, \tau_1)\overline{Z_{j}(v, \tau_2)} $ & $A_{k-1}$ \\
\hline

k=4n & $\sum\limits_{\substack{0\leq j<2n,\\ j:even}}(Z_j(u,\tau_1)+Z_{4n-j}(u,\tau_1))\overline{(Z_j(v,\tau_2)+Z_{4n-j}(v,\tau_2))}+2Z_{2n}(u, \tau_1)\overline{Z_{2n}(v, \tau_2)}$ & $D_{2n+2}$ \\
\hline
k=4n-2 & $\sum\limits_{\substack{0\leq j\leq 4n-2,\\ j:even}}Z_j(u,\tau_1)\overline{Z_j(v,\tau_2)}+Z_{2n-1}(u, \tau_1)\overline{Z_{2n-1}(v, \tau_2)}$ & $D_{2n+1}$ \\
&$+\sum\limits_{\substack{1\leq j\leq 2n-3,\\ j:odd}}\left(Z_j(u,\tau_1)\overline{Z_{4k-2-j}(v,\tau_2)}+\overline{Z_j(v,\tau_1)}{Z_{4k-2-j}(u,\tau_2)}\right)$&\\
\hline
k=10 & $(Z_0(u,\tau_1)+Z_6(u,\tau_1))\overline{(Z_0(v,\tau_2)+Z_6(v,\tau_2))}+(Z_6(u,\tau_1)+Z_{11}(u,\tau_1))\overline{(Z_6(v,\tau_2)+Z_{11}(v,\tau_2))}$ & $ E_6$ \\
&$+(Z_4(u,\tau_1)+Z_{10}(u,\tau_1))\overline{(Z_4(v,\tau_2)+Z_{10}(v,\tau_2))}$&\\
\hline
k=16 & $(Z_0(u,\tau_1)+Z_{16}(u,\tau_1))\overline{(Z_0(v,\tau_2)+Z_{16}(v,\tau_2))}+(Z_3(u,\tau_1)+Z_7(u,\tau_1))\overline{(Z_3(v,\tau_2)+Z_7(v,\tau_2))}$ & $ E_7$ \\
&$+(Z_4(u,\tau_1)+Z_{8}(u,\tau_1)+Z_{12}(u,\tau_1))\overline{(Z_4(v,\tau_2)+Z_{8}(v,\tau_2)+Z_{12}(v,\tau_1))}$&\\
\hline
 k=28 & $(Z_0(u,\tau_1)+Z_{10}(u,\tau_1)+Z_{18}(u,\tau_1)+Z_{28}(u,\tau_1))\overline{(Z_0(v,\tau_2)+Z_{10}(v,\tau_2)+Z_{18}(v,\tau_2)+Z_{28}(v,\tau_2))}$ & $ E_8$\\
 &$+(Z_6(u,\tau_1)+Z_{12}(u,\tau_1)+Z_{16}(u,\tau_1)+Z_{22}(u,\tau_1))\overline{(Z_6(v,\tau_2)+Z_{12}(v,\tau_2)+Z_{16}(v,\tau_2)+Z_{22}(v,\tau_2))}$&\\
 \hline
\end{tabular}
$$}
\end{theorem}

From Theorem \ref{invariant}, we know that  there is a modular invariant associated to an extension vertex operator algebra of $L_{{sl_2}}(k,0)$. In the following we will construct extension vertex operator algebras of $L_{{sl_2}}(k,0)$ which satisfy conditions (i) and (ii) and realize  modular invariants in Theorem \ref{CIZ1}.
\subsection{Realization of modular invariants of type $D$}In this subsection we  construct extension vertex operator algebras of $L_{{sl_2}}(k,0)$ which realize the modular invariants of type $D$. First, we recall some facts about simple current extension. Recall that if $V$ is a vertex operator algebra satisfying conditions (i) and (ii), then for any $V$-modules $W^1$ and $W^2$ the tensor product $W^1\boxtimes W^2$of $W^1$ and $W^2$ exists. An irreducible $V$-module $M$ is called a {\em simple current} $V$-module if for any irreducible $V$-module $W$, the tensor product $W\boxtimes M$ of $W$ and $M$ is irreducible. An extension vertex operator algebra $U$ of $V$ is called a {\em simple current extension} if $U$ is isomorphic to $V\oplus M^1\oplus ...\oplus M^k$ as $V$-module and   $M^i,\, 1\leq i\leq k$, are simple current $V$-modules.
\begin{theorem}\label{simple current}\cite{DM1}
Let $U^1$, $U^2$ be simple current extension vertex operator algebras of $V$. Assume that $U^1$ is isomorphic to  $U^2$ as $V$-module, then $U^1$ is isomorphic to $U^2$ as vertex operator algebra.
\end{theorem}

Recall from \cite{Li} that  $ L_{{sl_2}}(k,4n)$ is a simple current module of $L_{{sl_2}}(k,0)$ if $k=4n$ for some positive integer $n$. Furthermore, we have the following results  proved in \cite{Li}.
 \begin{theorem}\label{A}
If $k=4n$ for some positive integer $n$.

(1) $L_{{sl_2}}(k,0)\oplus L_{{sl_2}}(k,4n)$ has a vertex operator algebra structure. Moreover, this vertex operator algebra satisfies  conditions (i) and (ii).

(2) The vector spaces $L_{{sl_2}}(k,j)\oplus L_{{sl_2}}(k, k-j)$ for $0\leq j < 2n$ and $j$ is even,  have irreducible module structures for  $L_{{sl_2}}(k,0)\oplus L_{{sl_2}}(k,4n)$.

(3) The vector space $L_{{sl_2}}(k, 2n)$ has two nonisomorphic irreducible module structures for  $L_{{sl_2}}(k,0)\oplus L_{{sl_2}}(k,4n)$, which are denoted by $L_{{sl_2}}(k, 2n)$ and $L_{{sl_2}}(k, 2n)^{\sigma}$ respectively.

(4) Any irreducible module for $L_{{sl_2}}(k,0)\oplus L_{{sl_2}}(k,4n)$ is isomorphic to one of the modules in (2) and (3).

 \end{theorem}

Therefore, by  Theorem \ref{invariant}, we know that $L_{{sl_2}}(k,0)\oplus L_{{sl_2}}(k,4n)$ realizes a type $D$ modular invariant of  $L_{{sl_2}}(k, 0)$:
 \begin{tiny}
 \begin{align*}\sum\limits_{\substack{0\leq j<2n,\\ j:even}}(Z_j(u,\tau_1)+Z_{4n-j}(u,\tau_1))\overline{(Z_j(v,\tau_2)+Z_{4n-j}(v,\tau_2))}+2Z_{2n}(u, \tau_1)\overline{Z_{2n}(v, \tau_2)}.
 \end{align*}
 \end{tiny}
\subsection{Realization of modular invariants of type $E_6$}
In this subsection we  construct an extension vertex operator algebra of $L_{{sl_2}}(k,0)$ which realizes the modular invariant of type $E_6$. Consider the conformal inclusions $$SU(2)_{10}\subset Spin(5)_1\ \ \ \ \ \ \  (L_{sl_2}(10,0)\subset L_{B_2}(1,0)),$$
from the definition of conformal inclusion we know that $ L_{B_2}(1,0)$ and $L_{sl_2}(10,0)$ have the same conformal vector \cite{SW}, thus $ L_{B_2}(1,0)$ is an extension vertex operator algebra of $L_{sl_2}(10,0)$. The decomposition of $ L_{B_2}(1,0)$ as $L_{sl_2}(10,0)$-module is as follows:$$L_{B_2}(1,0)=L_{sl_2}(10, 0)\oplus L_{sl_2}(10, 6).$$ Since $ L_{B_2}(1,0)$ satisfies conditions (i) and (ii), from Theorem \ref{invariant} there exists a modular invariant of $L_{sl_2}(10,0)$ associated to $ L_{B_2}(1,0)$. Obviously, this modular invariant is equal to the  modular invariant of  type $E_6$:
\begin{tiny}
\begin{align*}
(Z_0&(u,\tau_1)+Z_6(u,\tau_1))\overline{(Z_0(v,\tau_2)+Z_6(v,\tau_2))}+(Z_3(u,\tau_1)+Z_7(u,\tau_1))\overline{(Z_3(v,\tau_2)+Z_7(v,\tau_2))}\\
&+(Z_4(u,\tau_1)+Z_{10}(u,\tau_1))\overline{(Z_4(v,\tau_2)+Z_{10}(v,\tau_2))}.
\end{align*}
\end{tiny}
Recall from \cite{K} that $L_{B_2}(1,0)$ has three irreducible modules $ L_{B_2}(\Lambda_0)$, $ L_{B_2}(\Lambda_1)$ and $ L_{B_2}(\Lambda_2)$, where $\Lambda_i$ denotes the fundamental weight of $\hat{B}_2$ and  $ L_{B_2}(\Lambda_i)$ denotes the highest weight module of  $\hat{B}_2$ with highest weight $\Lambda_i$. Thus from the modular invariant we have
\begin{align*}
&L_{B_2}(\Lambda_0)\cong L_{sl_2}(10, 0)\oplus L_{sl_2}(10, 6),\\
&L_{B_2}(\Lambda_1)\cong L_{sl_2}(10, 3)\oplus L_{sl_2}(10, 7),\\
&L_{B_2}(\Lambda_2)\cong L_{sl_2}(10, 4)\oplus L_{sl_2}(10, 10).
\end{align*}
Furthermore, it was proved in \cite{DJX} that the vertex operator algebra structure on $U=L_{sl_2}(10, 0)\oplus L_{sl_2}(10, 6)$ is unique.
\begin{theorem}\label{u1}
If $\overline U$ is an extension vertex operator algebra of $L_{sl_2}(10, 0)$ satisfying conditions (i) and (ii). Assume that $\overline U$ is isomorphic to $L_{sl_2}(10, 0)\oplus L_{sl_2}(10, 6)$ as $L_{sl_2}(10, 0)$-module, then $\overline U$ is isomorphic to $L_{B_2}(1, 0)$.
\end{theorem}
\subsection{Realization of modular invariants of type $E_8$}
In this subsection we  construct an extension vertex operator algebra of $L_{sl_2}(k,0)$ which realizes the modular invariant of type $E_8$. Consider the conformal inclusions $$SU(2)_{28}\subset (G_2)_1\ \ \ \ \ \ \ \ L_{sl_2}(28, 0)\subset L_{G_2}(1,0).$$ From the definition of conformal inclusion we know that $ L_{G_2}(1,0)$ and $L_{sl_2}(28,0)$ have the same conformal vector \cite{SW}, thus $ L_{G_2}(1,0)$ is an extension vertex operator algebra of $L_{sl_2}(28,0)$. The decomposition of $ L_{G_2}(1,0)$ as $L_{sl_2}(28,0)$-module is as follows:$$L_{G_2}(1,0)=L_{sl_2}(28, 0)\oplus L_{sl_2}(28, 10)\oplus L_{sl_2}(28, 18)\oplus L_{sl_2}(28, 28).$$Since $ L_{G_2}(1,0)$ satisfies conditions (i) and (ii), from Theorem \ref{invariant} there exists a modular invariant of $L_{sl_2}(28,0)$ associated to $ L_{G_2}(1,0)$. Obviously, this modular invariant is equal to the modular invariant of type $E_8$:
\begin{tiny}
\begin{align*}
 &(Z_0(u,\tau_1)+Z_{10}(u,\tau_1)+Z_{18}(u,\tau_1)+Z_{28}(u,\tau_1))\overline{(Z_0(v,\tau_2)+Z_{10}(v,\tau_2)+Z_{18}(v,\tau_2)+Z_{28}(v,\tau_2))}\\ &+(Z_6(u,\tau_1)+Z_{12}(u,\tau_1)+Z_{16}(u,\tau_1)+Z_{22}(u,\tau_1))\overline{(Z_6(v,\tau_2)+Z_{12}(v,\tau_2)+Z_{16}(v,\tau_2)+Z_{22}(v,\tau_2))}.
\end{align*}
\end{tiny}
Recall from \cite{K} that $L_{G_2}(1,0)$ has two irreducible modules $ L_{G_2}(\Lambda_0)$, $ L_{G_2}(\Lambda_2)$, where $\Lambda_i$ denotes the fundamental weight of $\hat{G}_2$ and  $ L_{G_2}(\Lambda_i)$ denotes the highest weight module of  $\hat{G}_2$ with highest weight $\Lambda_i$. Thus from the modular invariant we have\begin{align*}
&L_{G_2}(\Lambda_0)\cong L_{sl_2}(28, 0)\oplus L_{sl_2}(28, 10)\oplus L_{sl_2}(28, 18)\oplus L_{sl_2}(28, 28),\\
&L_{G_2}(\Lambda_2)\cong L_{sl_2}(28, 6)\oplus L_{sl_2}(28, 12)\oplus L_{sl_2}(28, 16)\oplus L_{sl_2}(28, 22).
\end{align*}

We now prove  that the vertex operator algebra structure on $U=L_{sl_2}(28, 0)\oplus L_{sl_2}(28, 10)\oplus L_{sl_2}(28, 18)\oplus L_{sl_2}(28, 28)$ is unique.
\begin{theorem}\label{u2}
If $\overline U$ is an extension vertex operator algebra of $L_{sl_2}(28, 0)$ satisfying conditions (i) and (ii). Assume that $\overline U$ is isomorphic to $L_{sl_2}(28, 0)\oplus L_{sl_2}(28, 10)\oplus L_{sl_2}(28, 18)\oplus L_{sl_2}(28, 28)$ as $L_{sl_2}(28, 0)$-module, then $\overline U$ is isomorphic to $L_{G_2}(1, 0)$.
\end{theorem}
\pf From the assumption, $\overline U$ is isomorphic to $L_{sl_2}(28, 0)\oplus L_{sl_2}(28, 10)\oplus L_{sl_2}(28, 18)\oplus L_{sl_2}(28, 28)$ as $L_{sl_2}(28, 0)$-module, thus the weight one subspace $\overline U_1$ of $\overline U$ is nonzero. Since $\overline U$ is a vertex operator algebra satisfying conditions (i) and (ii),  $\overline U_1$  is a reductive Lie algebra under the bracket operation $[a, b]=a_0b$ (see \cite{DM1}). By direct computation we immediately know that $sl_2$ is a Lie subalgebra of $\overline{U}_1$  and $\overline{U}_1$ is isomorphic to $sl_2\oplus L_{10}$ as $sl_2$-module, where $L_{10}$ is the highest weight module of $sl_2$ with highest weight $10$. Since $\overline{U}$ is self-dual, the restriction of the  commutator map $[,]:L_{10}\otimes L_{10}\rightarrow sl_2$ is nonzero, then  the Lie algebra $\overline{U}_1$ is isomorphic to $G_2$ by Lemma 6.7 in \cite{KO}. Thus we have $\overline U$ is a $\hat{G}_2$-module (see \cite{DM2}). Assume that  the level of $\overline U$ is equal to $k$, then  the vertex operator subalgebra $\langle \overline U_1\rangle$ generated by $ \overline U_1$ is isomorphic to the affine vertex operator algebra $L_{G_2}(k, 0)$ (see \cite{DM2}). Since $L_{G_2}(k, 0)$ and $\overline U$ have the same central charge,  $k$ must be equal to $1$, implies that $\langle \overline U_1\rangle$ is isomorphic to $L_{G_2}(1, 0)$. Note that $L_{G_2}(1, 0)$ is isomorphic to $L_{sl_2}(28, 0)\oplus L_{sl_2}(28, 10)\oplus L_{sl_2}(28, 18)\oplus L_{sl_2}(28, 28)$ as $L_{sl_2}(28, 0)$-module, this implies that $\overline{U}=\langle \overline U_1\rangle$ and then  the vertex operator algebra  $\overline{U}$ is isomorphic to $L_{G_2}(1, 0)$. The proof is complete.
\qed
\subsection{Classification of extension vertex operator algebras of $L_{sl_2}(k,0)$}

In this subsection, we classify extension vertex operator algebras of $L_{sl_2}(k,0),$
which is the first main result in this paper.
\begin{theorem}\label{kvoa13}
Let $U$ be an extension  vertex operator algebra of $L_{sl_2}(k,0)$ which satisfies  conditions (i) and (ii). Then $U$ is isomorphic to either $L_{sl_2}(k,0)$ or one of vertex operator algebras in Theorems \ref{A},  \ref{u1}, \ref{u2}.
\end{theorem}
\pf From Theorem \ref{invariant},  $U$  realizes a modular invariant of $L_{sl_2}(k,0)$. We first prove that this modular invariant is not equal to the  modular invariant of type $D_{2n+1}$.  Otherwise, if the modular invariant is equal to the modular invariant of type $D_{2n+1}$, it follows that $U$ should be isomorphic to $L_{sl_2}(k,0)$, but the the modular invariant of type $D_{2n+1}$ is not equal to the modular invariant of type $A_{k-1}$, which is the modular invariant associated to the extension vertex operator algebra $L_{sl_2}(k,0)$, this is a contradiction. We now prove that the modular invariant is not equal to the modular invariant of type  $E_7$. Otherwise, if the modular invariant is  equal to the modular invariant of type  $E_7$, then we have $U$ is isomorphic to $L_{sl_2}(16, 0)\oplus L_{sl_2}(16, 16)$ as $L_{sl_2}(16, 0)$-module. Note that $L_{sl_2}(16, 16)$ is a simple current module of $L_{sl_2}(16, 0)$, from Theorem \ref{simple current} we know that the vertex operator algebra structure on $L_{sl_2}(16, 0)\oplus L_{sl_2}(16, 16)$ is unique, this implies that the modular invariant is also equal to the modular invariant of type $D_{10}$, this is a contradiction. So the modular invariant is equal to one of the other modular invariants in  Theorem \ref{CIZ1}. By Theorems \ref{simple current}, \ref{u1},  \ref{u2}, $U$ is  isomorphic to either $L_{sl_2}(k,0)$ or one of vertex operator algebras in Theorems \ref{A},  \ref{u1}, \ref{u2}. The proof is complete.  \qed
\begin{remark}
The rigid associative commutative algebras in the category of $U_q(sl_2)$-modules were classified in \cite{KO}. Using a recent result in \cite{HKL}, Theorem \ref{kvoa13} is essentially obtained in \cite{KO}.
\end{remark}

\begin{remark}
From Theorem \ref{kvoa13}, we know that the modular invariants of type $D_{2n+1}$ can not be realized by extension vertex operator algebras of $L_{sl_2}(4n-2,0)$. However, the modular invariants of type $D_{2n+1}$ were realized in \cite{DH} by extension vertex operator superalgebras of $L_{sl_2}(4n-2,0)$.
\end{remark}
\section{Classification of preunitary vertex operator algebras with central charge $c<1$}
\def\theequation{4.\arabic{equation}}
\setcounter{equation}{0}

In this section we  classify  preunitary  vertex operator algebras with central charge $c<1$.
 First, we recall some facts about
Virasoro vertex operator algebras \cite{FZ}, \cite{W}. We denote the Virasoro algebra by  $L=\oplus_{n\in \mathbb{Z}}
\mathbb{C}L_n\oplus \mathbb{C}C$ with the
commutation relations
\begin{align*}
[L_m, L_n]=(m-n)&L_{m+n}+\frac{1}{12}(m^3-m)\delta_{m+n,
0}C,\\
&[L_m, C]=0.
\end{align*}

Set $\mathfrak{b}=(\oplus_{n\geq
1}\mathbb{C}L_n)\oplus(\mathbb{C}L_0\oplus \mathbb{C}C)$, then  $\mathfrak{b}$ is a subalgebra of $L$. For any two
complex numbers $c, h\in \C$, let $\C$ be a 1-dimensional
$\mathfrak{b}$-module defined as follows:
\begin{align*}
&L_n\cdot 1=0, n\geq 1,\\
&L_0\cdot 1=h\cdot 1,\\
&C\cdot 1=c\cdot 1.
\end{align*}
Set $$V(c, h)=U(L)\otimes_{U(\mathfrak{b})}\C$$ where
$U(\cdot)$ denotes the universal enveloping algebra. Then $V(c, h)$ is
a highest weight module of the Virasoro algebra of highest weight $(c, h)$, which is called the Verma module of Virasoro algebra,
and $V(c, h)$ has a unique maximal proper submodule $J(c, h)$. Let
$L(c, h)$ be the unique irreducible quotient module of $V(c, h)$.
Set
$$\overline{V(c, 0)}=V(c, 0)/(U(L)L_{-1}1\otimes 1),$$ it is well-known
that $\overline{V(c, 0)}$ has a vertex operator algebra structure with Virasoro element $\w=L_{-2}1$
and $ L(c, 0)$ is the unique irreducible quotient vertex operator
algebra of $\overline{V(c, 0)}$ \cite{FZ}.

For $m\geq 0$, set
$$c_m=1-\frac{6}{(m+2)(m+3)},$$ $$h^m_{r,
s}=\frac{(r(m+3)-s(m+2))^2-1}{4(m+2)(m+3)},\  (1\leq s\leq r \leq m+1).$$
It was proved in \cite{W}, \cite{DLM2} that $L(c_m, 0)$ $(m\geq 0)$ are rational and
$C_2$-cofinite vertex operator algebras and $L(c_m, h^m_{r, s})$ are
the complete list of irreducible $L(c_m, 0)$-modules.
In \cite{CIZ1}, \cite{CIZ2} and \cite{Ka},  modular invariants of Virasoro vertex operator algebra $L(c_m, 0)$ were classified:
\begin{theorem}\label{ciz2}
Let $m$ be a nonnegative integer. We use $Z_{r,s}(u,\tau)$ to denote the trace function associated to $L(c_m, 0)$-module $L(c_m, h^m_{r,s})$. Then any modular invariant of Virasoro vertex operator algebra $L(c_m, 0)$ is equal to  one of the following modular invariants:\\
\begin{tiny}
\begin{tabular}{|c c c|}
\hline

any m  & $\sum\limits_{1\leq s\leq r\leq m-1}Z_{r,s}(u,\tau_1)\overline{Z_{r,s}(v,\tau_2)}$ & $(A_{m-1}, A_{m})$ \\
\hline
m=4n & $\frac{1}{2}\bigg\{\sum\limits_{\substack{1\leq s\leq r \leq 4n+1,\\r:odd}}(Z_{r,s}(u,\tau_1)+Z_{4n+2-r, s}(u,\tau_1))\overline{(Z_{r,s}(v,\tau_2)+Z_{4n+2-r, s}(v,\tau_2))}\bigg\}$&$(D_{2n+2}, A_{4n+2})$ \\
\hline
m=4n-2 & $\frac{1}{2}\sum\limits_{s=1}^{4n}\bigg\{\sum\limits_{r:odd}Z_{r,s}(u,\tau_1)\overline{Z_{r,s}(v,\tau_2)}+Z_{2n, s}(u, \tau_1)\overline{Z_{2n, s}(v, \tau_2)}$&$(D_{2n+1}, A_{4n})$ \\
& $+\sum\limits_{r:even}Z_{r, s}(u,\tau_1)\overline{Z_{4n-r, s}(v,\tau_2)}\bigg\} $ & \\
\hline
m=4n-1 & $\frac{1}{2}\bigg\{\sum\limits_{\substack{1\leq s\leq r \leq 4n,\\s:odd}}(Z_{r,s}(u,\tau_1)+Z_{r, 4n+2-s}(u,\tau_1))\overline{(Z_{r,s}(v,\tau_2)+Z_{r,4n+2- s}(v,\tau_2))}\bigg\}$ & $(A_{4n}, D_{2n+2})$ \\
\hline
m=4n-3 & $\frac{1}{2}\sum\limits_{s=1}^{4n-2}\bigg\{\sum\limits_{s:odd}Z_{r,s}(u,\tau_1)\overline{Z_{r,s}(v,\tau_2)}+Z_{r, 2n}(u, \tau_1)\overline{Z_{r, 2n}(v, \tau_2)}$&$( A_{4n-2},D_{2n+1})$ \\
& $+\sum\limits_{s:even}Z_{r, s}(u,\tau_1)\overline{Z_{r, 4n-s}(v,\tau_2)}\bigg\} $ & \\
\hline
m=9 & $\frac{1}{2}\sum\limits_{r=1}^{10}\bigg\{ (Z_{r,1}(u,\tau_1)+Z_{r,7}(u,\tau_1))\overline{(Z_{r,1}(v,\tau_2)+Z_{r,7}(v,\tau_2))}$&$(A_{10}, E_6)$ \\
&$+(Z_{r,4}(u,\tau_1)+Z_{r,8}(u,\tau_1))\overline{(Z_{r,4}(v,\tau_2)+Z_{r,8}(v,\tau_2))}$&\\
&$+(Z_{r,5}(u,\tau_1)+Z_{r,11}(u,\tau_1))\overline{(Z_{r,5}(v,\tau_2)+Z_{r,11}(v,\tau_2))}\bigg\}$ & \\
\hline
m=10 & $\frac{1}{2}\sum\limits_{s=1}^{12}\bigg\{ (Z_{1, s}(u,\tau_1)+Z_{7, s}(u,\tau_1))\overline{(Z_{1, s}(v,\tau_2)+Z_{7, s}(v,\tau_2))}$&$(E_6, A_{12})$\\
&$+(Z_{4, s}(u,\tau_1)+Z_{8, s}(u,\tau_1))\overline{(Z_{4, s}(v,\tau_2)+Z_{8, s}(v,\tau_2))}$&\\
&$+(Z_{5, s}(u,\tau_1)+Z_{11, s}(v,\tau_1))\overline{(Z_{5, s}(v,\tau_2)+Z_{11, s}(v,\tau_2))}\bigg\}$ & \\
\hline
m=15 & $\frac{1}{2}\sum\limits_{r=1}^{16}\bigg\{ (Z_{r, 1}(u,\tau_1)+Z_{r,17}(u,\tau_1))\overline{(Z_{r, 1}(v,\tau_2)+Z_{r,17}(v,\tau_2))}$&$(A_{16}, E_{7})$\\
&$+(Z_{r, 5}(u,\tau_1)+Z_{r, 13}(u,\tau_1))\overline{(Z_{r, 5}(v,\tau_2)+Z_{r, 13}(v,\tau_2))}$&\\
&$+(Z_{r, 7}(u,\tau_1)+Z_{r, 11}(u,\tau_1))\overline{(Z_{r, 7}(v,\tau_2)+Z_{r, 11}(v,\tau_2))}$ & \\
&$+(Z_{r, 3}(u,\tau_1)+Z_{r,15}(u,\tau_1))\overline{Z_{r,9}(v,\tau_2)}+Z_{r,9}(v,\tau_2)\overline{(Z_{r, 3}(u,\tau_1)+Z_{r,15}(u,\tau_1))}\bigg\}$&\\
\hline
m=16 &  $\frac{1}{2}\sum\limits_{s=1}^{18}\bigg\{ (Z_{1, s}(u,\tau_1)+Z_{17,s}(u,\tau_1))\overline{(Z_{1, s}(v,\tau_2)+Z_{17,s}(v,\tau_2))}$&$(E_7, A_{18})$\\
&$+(Z_{ 5,s}(u,\tau_1)+Z_{ 13,s}(u,\tau_1))\overline{(Z_{ 5,s}(v,\tau_2)+Z_{ 13,s}(v,\tau_2))}$&\\
&$+(Z_{ 7,s}(u,\tau_1)+Z_{11,s}(u,\tau_1))\overline{(Z_{ 7,s}(v,\tau_2)+Z_{11,s}(v,\tau_2))}$ & \\
&$+(Z_{ 3,s}(u,\tau_1)+Z_{15,s}(u,\tau_1))\overline{Z_{9,s}(v,\tau_2)}+Z_{9,s}(v,\tau_2)\overline{(Z_{ 3,s}(u,\tau_1)+Z_{15,s}(u,\tau_1))}\bigg\}$&\\
\hline
 m=27 & $\frac{1}{2}\sum\limits_{r=1}^{28}\bigg\{(Z_{r,1}(u,\tau_1)+Z_{r,11}(u,\tau_1)+Z_{r,19}(u,\tau_1)+Z_{r,29}(u,\tau_1))$&$(A_{28}, E_8)$\\
 & $\overline{(Z_{r,1}(v,\tau_2)+Z_{r,11}(v,\tau_2)+Z_{r,19}(v,\tau_2)+Z_{r,29}(v,\tau_2))}$ & \\
 &$+(Z_{r,7}(u,\tau_1)+Z_{r,13}(u,\tau_1)+Z_{r,17}(u,\tau_1)+Z_{r,23}(u,\tau_1))$ &\\
 &$\overline{(Z_{r,7}(v,\tau_2)+Z_{r,13}(v,\tau_2)+Z_{r,17}(v,\tau_2)+Z_{r,23}(v,\tau_2))}\bigg\}$ &\\
 \hline
 m=28 & $\frac{1}{2}\sum\limits_{s=1}^{30}\bigg\{(Z_{1, s}(u,\tau_1)+Z_{11, s}(u,\tau_1)+Z_{19, s}(u,\tau_1)+Z_{29, s}(u,\tau_1))$&$( E_8, A_{30})$\\
 &$\overline{(Z_{1, s}(v,\tau_2)+Z_{11, s}(v,\tau_2)+Z_{19, s}(v,\tau_2)+Z_{29, s}(v,\tau_2))}$&\\
 &$+(Z_{7, s}(u,\tau_1)+Z_{13, s}(u,\tau_1)+Z_{17,s}(u,\tau_1)+Z_{23, s}(u,\tau_1))$&\\
 &$\overline{(Z_{7, s}(v,\tau_2)+Z_{13, s}(v,\tau_2)+Z_{17, s}(v,\tau_2)+Z_{23, s}(v,\tau_2))}\bigg\}$&\\
 \hline
\end{tabular}
\end{tiny}
\end{theorem}

It is clear that vertex operator algebra $L(c_m,0)$ realizes the modular invariant of type $(A_{m-1}, A_{m}).$
So our discussion below only deals with the other types.

\subsection{Realization of modular invariants of type $(D_{2n+2}, A_{4n+2})$ and $(A_{4n}, D_{2n+2})$}
In this subsection we  construct extension vertex operator algebras of $L(c_m,0)$ which realize the modular invariants of type $(D_{(m/2)+1}, A_m)$ and $(A_{m-1}, D_{(m+3)/2})$.  Let $A_1=\Z\alpha$, with $\langle\alpha, \alpha\rangle=2$, be the root lattice of type $A_1$. If $m$ is a positive integer such that  $m\equiv 0$ mod $4.$ Let $A^{m+1}_
1 = \Z\alpha^0\oplus \Z\alpha^1\oplus\cdots\oplus \Z\alpha^m$ be the
orthogonal sum of $m+1$ copies of $A_1$ and set $L=A_1^{m+1}\cup(\gamma_1+A_1^{m+1})$ where $\gamma_1=(1/2)
\sum^{m-1}_{
i=0} \alpha^i.$ Consider the lattice vertex operator algebra $V_L$ associated to $L$. It was proved in \cite{LLY} that $L(c_1, 0)
\otimes\cdots \otimes L(c_{m-1}, 0)\otimes L_{sl_{2}}(m+ 1, 0)$  is a vertex operator subalgebra of $V_L$ and the commutant vertex operator algebra of $L(c_1, 0)
\otimes\cdots \otimes L(c_{m-1}, 0)\otimes L_{sl_{2}}(m+ 1, 0)$ in $V_L$ is isomorphic to $L(c_m, 0)\oplus L(c_m, h_{1, m+2}^m)$ as $L(c_m, 0)$-module. Furthermore,  the following results were obtained in \cite{LLY}.
\begin{theorem}\label{vvoa1}
  Let $m$ be a positive integer such that  $m\equiv 0 \mod 4.$

(1) There is a  vertex operator algebra structure on  $U=L(c_m, 0)\oplus L(c_m, h_{1, m+2}^m)$. Moreover, $U$ satisfies conditions (i) and (ii).

(2) For any  $1\leq r<(m+2)/2, 1\leq s<(m+3)/2$ and $r$ is odd, $L(c_m, h_{r,s}^m)\oplus L(c_m, \tilde{h}_{r,s}^m)$ is an irreducible $U$-module, where $\tilde{h}_{r,s}^m=h_{m+2-r, s}^m$.

(3) For $1\leq s<(m+3)/2,$ $L(c_m, h_{(m+2)/2, s})$ has exactly two nonisomorphic irreducible $U$-module structures denoted by $L(c_m, h_{(m+2)/2, s})$ and $L(c_m, h_{(m+2)/2, s})^{\sigma}$ respectively.

(4) Any irreducible $U$-module is isomorphic to one of irreducible $U$-modules in (2), (3).
\end{theorem}
Thus, by Theorem \ref{invariant}, $U=L(c_m, 0)\oplus L(c_m, h_{1, m+2}^m)$ 
realizes a  modular invariant of  $L(c_m, 0)$:
\begin{tiny}
\begin{align*}
\sum_{\substack{1\leq r<(m+2)/2,\\ 1\leq s<(m+3)/2,\\r:odd}}(Z_{{r,s}}(u, \tau_1)+Z_{{m+2-r,s}}(u, \tau_1))\overline{(Z_{{r,s}}(v, \tau_2)+Z_{{m+2-r,s}}(u, \tau_2))}+2Z_{{\frac{m+2}{2},s}}(u, \tau_1)\overline{Z_{{\frac{m+2}{2},s}}(v, \tau_2)}.
\end{align*}
\end{tiny}
Note that $h_{r, s}=h_{m+2-r, m+3-s}$, then this modular invariant is equal to the  modular invariant of type $(D_{\frac{m}{2}+2}, A_{m+2})$. Moreover, we have
\begin{theorem}\label{vvoa3}
The vertex operator algebra $U=L(c_m, 0)\oplus L(c_m, h_{1, m+2}^m)$ in Theorem \ref{vvoa1} is a unitary vertex operator algebra.
\end{theorem}
\pf We are working in the setting of \cite{DLin}, \cite{LLY}. It is known  \cite{DLin} that $(V_L, \phi)$ is a unitary vertex operator algebra, where $\phi$ is an anti-linear involution of $V_L$. To prove $U=L(c_m, 0)\oplus L(c_m, h_{1, m+2}^m)$ is a unitary vertex operator algebra, by Corollary 2.8 in \cite{DLin} we only need to prove that $\phi(\w')=\w'$, where $\w'$ denotes the conformal vector of $L(c_1, 0)
\otimes\cdots \otimes L(c_{m-1}, 0)\otimes L_{sl_{2}}(m+ 1, 0)$. This follows immediately from the definition of $\phi$ in \cite{DLin} and formula (3.4) in \cite{LLY}. This finishes the proof.
\qed
\vskip0.5cm
We now assume $m$ is a positive integer such that  $m\equiv 3\mod 4.$ 
Let $L=A_1^{m+1}\cup(\gamma_2+A_1^{m+1})$ where $\gamma_2=(1/2)
\sum^m_{
i=0} \alpha^i.$ Consider the lattice vertex operator algebra $V_L$ associated to $L$. It was proved in \cite{LLY} that $L(c_1, 0)
\otimes\cdots \otimes L(c_{m-1}, 0)\otimes L_{sl_{2}}(m+ 1, 0)$  is a vertex operator subalgebra of $V_L$ and the commutant vertex operator algebra of $L(c_1, 0)
\otimes\cdots \otimes L(c_{m-1}, 0)\otimes L_{sl_{2}}(m+ 1, 0)$ in $V_L$ is isomorphic to $L(c_m, 0)\oplus L(c_m, h_{1, m+2}^m)$ as $L(c_m, 0)$-module. Furthermore, the following results  were obtained in \cite{LLY}.
\begin{theorem}\label{vvoa2}
Let $m$ be a positive integer such that  $m\equiv 3\mod 4.$ Then\\
(1) There is a  vertex operator algebra structure on  $U=L(c_m, 0)\oplus L(c_m, h_{1, m+2}^m)$. Moreover, $U$ satisfies conditions (i) and (ii).\\
(2) For any $1\leq r<(m+2)/2, 1\leq s<(m+3)/2$ and $s$ is odd, $L(c_m, h_{r,s}^m)\oplus L(c_m, \tilde{h}_{r,s}^m)$ is an irreducible $U$-module,  where $\tilde{h}_{r,s}^m=h_{r, m+3-s}^m$.\\
(3) For $1\leq s<(m+3)/2,$ $L(c_m, h_{r, (m+3)/2})$ has exactly two nonisomorphic irreducible $U$-module structures denoted by $L(c_m, h_{r,(m+3)/2})$ and $L(c_m, h_{r, (m+3)/2})^{\sigma}$ respectively.\\
(4) Any irreducible $U$-module is isomorphic to one of irreducible $U$-modules in (2), (3).
\end{theorem}

Thus, by Theorem \ref{invariant},  $U=L(c_m, 0)\oplus L(c_m, h_{1, m+2}^m)$ 
realizes a  modular invariant of  $L(c_m, 0)$:
\begin{tiny}
\begin{align*}
\sum_{\substack{1\leq r<(m+2)/2,\\ 1\leq s<(m+3)/2,\\s:odd}}(Z_{{r,s}}(u, \tau_1)+Z_{{r,m+3-s}}(u, \tau_1))\overline{(Z_{{r,s}}(v, \tau_2)+Z_{{r,m+3-s}}(v, \tau_2))}+2Z_{{r, \frac{m+3}{2}}}(u, \tau_1)\overline{Z_{{r, \frac{m+3}{2}}}(v, \tau_2)}.
\end{align*}
\end{tiny}
By direct computation, this modular invariant is equal to the  modular invariant of  type $(A_{m+1}, D_{\frac{m+3}{2}+1})$.
Moreover, by an argument similar to  the proof of Theorem \ref{vvoa3}, we have
\begin{theorem}
The vertex operator algebra $U=L(c_m, 0)\oplus L(c_m, h_{1, m+2}^m)$ in Theorem \ref{vvoa2} is a unitary vertex operator algebra.
\end{theorem}

\subsection{Realization of modular invariants of type $( E_6, A_{12})$ and $( E_8, A_{30})$}
In this subsection we  construct extension vertex operator algebras of $L(c_m,0)$ which realize the modular invariants of type $(A_{10}, E_6)$ and $(A_{28}, E_8)$.  First,  we recall some facts about the coset construction. Let $\{h, e, f \}$ be a standard basis of $sl_2$. For any irreducible $L_{sl_2}(m, 0)$-module $L_{sl_{2}}(m, j)$, consider the Sugawara operators $$L_{L_{sl_2}(m, 0)}(k)=\frac{1}{2(m+2)}\sum\limits_{j\in\Z}(:h_{-j}h_{k+j}:+:e_{-j}f_{k+j}:+:f_{-j}e_{k+j}:),$$ where $:\ \ :$ denotes the normal ordered product. These operators define a natural Virasoro action on $L_{sl_{2}}(k, j)$ with central charge $3m/(m+2)$.

Let $L_{sl_{2}}(m,i)$ be an irreducible $L_{sl_2}(m, 0)$-module and $L_{sl_{2}}(1,j)$ be an irreducible $L_{sl_2}(1, 0)$-module. Then $\hat{sl}_2$ acts on $L_{sl_{2}}(m,i)\otimes L_{sl_{2}}(1,j)$
by the diagonal action
$$x(n)\cdot (v \otimes w) =x_nv\otimes w + v \otimes x_nw$$
for any $x(n)\in \hat{sl}_2$ and $v \otimes w \in L_{sl_{2}}(m,i)\otimes L_{sl_{2}}(1,j)$. This gives a level $m + 1$ representation of $\hat{sl}_2$, thus $L_{sl_{2}}(m,i)\otimes L_{sl_{2}}(1,j)$ has a $L_{sl_2}(m+1, 0)$-module structure. The Sugawara operators $L_{L_{sl_2}(m+1, 0)}(
k)$ form a Virasoro algebra on $L_{sl_{2}}(m,i)\otimes L_{sl_{2}}(1,j)$ with central
charge $3(m+ 1)/(m+ 3)$. Set $L(k) = L_{L_{sl_2}(m+1, 0)}(k)-  L_{L_{sl_2}(m, 0)}(k)$, it is  known (cf. \cite{GKO}, \cite{KR}) that $L(k)$, $k \in  \Z$,
commute with the diagonal  operators $x(n)$ for all $x(n)\in \hat{sl}_2$ and they give rise to a
representation of the Virasoro algebra on $L_{sl_{2}}(m,i)\otimes L_{sl_{2}}(1,j)$ with central charge
$$c= c_m = 1 +\frac{3m}{(m+ 2)} -\frac{3(m+ 1)}{(m+ 3)} = 1 -\frac{6}{(m+ 2)(m+ 3)}.$$
Thus  $L_{sl_{2}}(m,i)\otimes L_{sl_{2}}(1,j)$ has a $L_{sl_2}(m+1, 0)\otimes L(c_m, 0)$-module structure. Since $L_{sl_2}(m+1, 0)\otimes L(c_m, 0)$ is rational, $L_{sl_{2}}(m,i)\otimes L_{sl_{2}}(1,j)$ is completely reducible as $L_{sl_2}(m+1, 0)\otimes L(c_m, 0)$-module, the decomposition is as follows:
\begin{align*}
L_{sl_{2}}(m,n)\otimes  L_{sl_{2}}(1, \epsilon ) &=\mathop\oplus\limits_{\substack{0\leq s\leq n,\\  s\equiv n+\epsilon \mod 2}}L(c_m, h^m_{n+1,s+1})\otimes L_{sl_{2}}(m+ 1, s)\\
&\bigoplus\mathop\oplus\limits_{\substack{
n+1\leq s\leq m+1,\\
s\equiv n+\epsilon  \mod 2}}
L(c_m, h^m_{m-n+1,m+2-s})
\otimes  L_{sl_{2}}(m+ 1, s)\\
&=\mathop\oplus\limits_{\substack{
0\leq s\leq m+1,\\
s\equiv n+\epsilon \mod  2}}
L(c_m, h^m_{n+1,s+1})\otimes  L_{sl_{2}}(m+ 1, s)
\end{align*}
for any $\epsilon = 0, 1,$ and $0 \leq n \leq m.$ Note that $h^m_{ n+1,s+1} = h^m_{m-n+1,m+2-s}.$

\vskip0.5cm
We now begin to construct extension vertex operator algebras of $L(c_m,0)$  realizing the modular invariants of type $(A_{10}, E_6)$ and $(A_{28}, E_8)$. Consider the conformal inclusions $$SU(2)_{10}\subset Spin(5)_1\ \ \ \ \ \ \  (L_{sl_2}(10,0)\subset L_{B_2}(1,0)).$$ From the discussion in Section 3.2, we know that the decomposition of $ L_{B_2}(1,0)$ as $L_{sl_2}(10,0)$-module is as follows:$$L_{B_2}(1,0)=L_{sl_2}(10, 0)\oplus L_{sl_2}(10, 6).$$
Consider the vertex operator algebra $L_{B_2}(1, 0) \otimes L_{sl_2}(1, 0)$, this vertex operator algebra contains $L_{sl_2}(10, 0)\otimes L_{sl_2}(1,0)$. The decomposition of $L_{B_2}(1, 0) \otimes L_{sl_2}(1, 0)$ as  $L_{sl_2}(10, 0)\otimes L_{sl_2}(1,0)$-module is as follows:
$$L_{B_2}(1,0) \otimes L_{sl_2}(1, 0)=L_{sl_2}(10, 0) \otimes L_{sl_2}(1, 0)\oplus L_{sl_2}(10, 6) \otimes L_{sl_2}(1, 0).$$
By the discussion above, we know that $L_{sl_2}(11, 0)$ is a vertex operator subalgebra of $L_{B_2}(1, 0) \otimes L_{sl_2}(1, 0)$, and we have
\begin{align*}
L_{B_2}(1,0) \otimes L_{sl_2}(1, 0)&=\mathop\oplus\limits_{\substack{
0\leq s\leq 11,\\
s\equiv 0 \mod  2}}
L(c_{10}, h^{10}_{1,s+1})\otimes  L_{sl_{2}}(11, s)\\
&\bigoplus\mathop\oplus\limits_{\substack{
0\leq s\leq 11,\\
s\equiv 0 \mod  2}}
L(c_{10}, h^{10}_{7,s+1})\otimes  L_{sl_{2}}(11, s).
\end{align*}
Thus the commutant vertex operator algebra of $L_{sl_2}(11, 0)$ in $L_{B_2}(1, 0) \otimes L_{sl_2}(1, 0)$ is $$U=L(c_{10}, 0)\oplus L(c_{10},h_{7,1}^{10}).$$
\begin{theorem}\label{vvoa5}
Let $U=L(c_{10}, 0)\oplus L(c_{10},h_{7,1}^{10})$ be the vertex operator algebra defined above. Then there is a unitary vertex operator algebra structure on $U$. Moreover, $U$ satisfies conditions (i) and (ii).
\end{theorem}
\pf  
  It is known \cite{DL} that $L_{B_2}(1,0) \otimes L_{sl_2}(1, 0)$ has a unitary vertex operator algebra structure $(L_{B_2}(1,0) \otimes L_{sl_2}(1, 0),\phi)$, where $\phi$ is an anti-linear involution of $L_{B_2}(1,0) \otimes L_{sl_2}(1, 0)$. By Corollary 2.8 in \cite{DLin}, to prove that the commutant vertex operator algebra of $L_{sl_2}(11, 0)$ in $L_{B_2}(1, 0) \otimes L_{sl_2}(1, 0)$ is a unitary vertex operator algebra,  we only need to prove that $\phi(\w')=\w'$, where $\w'$ denotes the conformal vector of $L_{sl_2}(11, 0)$. This follows immediately from the definitions of $\phi$ and the conformal vector of $L_{sl_2}(11, 0)$. Thus $U$ has a unitary vertex operator algebra structure. This implies that $U$ is a simple vertex operator algebra. In fact, note that  $\1\notin U^1$ if $U^1$ is a nontrivial ideal of $U$, then  for any $u\in U, u^1\in U^1$,
  \begin{align*}
  \langle u, u^1 \rangle&=\Res_zz^{-1}\langle Y_U(u,z)\1, u^1 \rangle\\
  &=\Res_zz^{-1}\langle \1, Y_U(e^{zL(1)}(-z^{-2})^{L(0)}\phi(u), z^{-1})u^1 \rangle\\
  &=0,
  \end{align*}
  this is a contradiction.

We now prove that $U$ is self dual. Note that $V'$ is a $V$-submodule of $U'$ and that $V$ is self dual, then $U'$ has a vacuum-like vector $\1$, this implies that $U'$ is isomorphic to $U$ as $U$-module, that is, $U$ is self dual. Therefore,   $U=L(c_{10}, 0)\oplus L(c_{10},h_{7,1}^{10})$ satisfies conditions (i) and (ii) by  Theorems \ref{C 2}, \ref{cvoa2}.
The proof is complete.
\qed
\vskip0.25cm
Therefore, by Theorem \ref{invariant},  $U=L(c_{10}, 0)\oplus L(c_{10},h_{7,1}^{10})$ realizes a   modular invariant of  $L(c_{10}, 0)$. Obviously, this modular invariant is equal to the modular invariant of  type $( E_6, A_{12})$:
\begin{tiny}
\begin{align*}
\frac{1}{2}\sum\limits_{p=1}^{12}&\bigg\{(Z_{1,q}(u,\tau_1)+Z_{7,q}(u,\tau_1))\overline{(Z_{1,q}(v,\tau_2)+Z_{7,q}(v,\tau_2))}+(Z_{4,q}(u,\tau_1)+Z_{8,q}(u,\tau_1))\overline{(Z_{4,q}(v,\tau_2)+Z_{8,q}(v,\tau_2))}\\
&+(Z_{5,q}(u,\tau_1)+Z_{11,q}(u,\tau_1))\overline{(Z_{5,q}(v,\tau_2)+Z_{11,q}(v,\tau_2))}\bigg\}.
\end{align*}
\end{tiny}
From the modular invariant, we immediately know that  the vector spaces
\begin{align*}
&L(c_{10}, h_{1,q}^{10})\oplus L(c_{10}, h_{7,q}^{10}),\\
&L(c_{10}, h_{4,q}^{10})\oplus L(c_{10}, h_{8,q}^{10}),\\
&L(c_{10}, h_{5,q}^{10})\oplus L(c_{10}, h_{11,q}^{10}),\ 1\leq q\leq 7,
 \end{align*}have irreducible $U$-module structures and these are exactly all the irreducible $U$-modules.

We now prove that the vertex operator algebra structure on $U=L(c_{10}, 0)\oplus L(c_{10},h_{7,1}^{10})$ is unique.
\begin{theorem}\label{vvoa4}
Let $\overline U$ be an extension vertex operator algebra of $L(c_{10}, 0)$ such that $\overline U$ is isomorphic to $L(c_{10}, 0)\oplus L(c_{10},h_{7,1}^{10})$ as $L(c_{10}, 0)$-module. Then $\overline U$ is isomorphic to the vertex operator algebra defined above.
\end{theorem}
\pf   Let  $M^1, M^2, M^3$  be irreducible $L_{sl_2}(10, 0)\otimes L_{sl_2}(1, 0)$-modules which are isomorphic to $L_{sl_2}(10, 0)\otimes L_{sl_2}(1, 0)$ or $L_{sl_2}(10, 6)\otimes L_{sl_2}(1, 0)$. From the discussion at the beginning of this subsection, we know that if $M^i$ is isomorphic to $L_{sl_2}(10, 0)\otimes L_{sl_2}(1, 0)$ (resp. $L_{sl_2}(10, 6)\otimes L_{sl_2}(1, 0)$) then $M^i$ has a $L(c_{10}, 0)$-submodule $\dot{M}^i$ which is isomorphic to $L(c_{10}, 0)$ (resp. $L(c_{10}, h_{7,1}^{10})$). Note that  $N_{M^1, M^2}^{M^3}\leq 1$,
 let $\Y_{M^1, M^2}^{M^3}$ be an intertwining operator  of type $\left(\begin{tabular}{c}
${M}^3$\\
${M}^1$ ${M}^2$\\
\end{tabular}\right)$ if $N_{M^1, M^2}^{M^3}=1$. Then  $\Y_{M^1, M^2}^{M^3}$ restricting to $\dot{M}^1, \dot{M}^2$ induces an intertwining operator $\Y_{\dot{M}^1, \dot{M}^2}^{\dot{M}^3}$ of type $\left(\begin{tabular}{c}
$\dot{M}^3$\\
$\dot{M}^1$ $\dot{M}^2$\\
\end{tabular}\right)$. Recall from \cite{W} that $N_{\dot{M}^1, \dot{M}^2}^{\dot{M}^3}\leq 1$, thus the vertex operators $\dot{Y}_{\overline U}$ of $\overline U$ is of the form:
$$\dot{Y}_{\overline U}|_{\dot{M}^1\times \dot{M}^2}=\lambda_{\dot{M}^1, \dot{M}^2}^{\dot{L}(c_{10}, 0)}\Y_{\dot{M}^1, \dot{M}^2}^{\dot{L}(c_{10}, 0)}+\lambda_{\dot{M}^1, \dot{M}^2}^{\dot{L}(c_{10}, 6)}\Y_{\dot{M}^1, \dot{M}^2}^{\dot{L}(c_{10}, 6)},$$
where $\lambda_{\dot{M}^1, \dot{M}^2}^{\dot{L}(c_{10}, 0)}$, $\lambda_{\dot{M}^1, \dot{M}^2}^{\dot{L}(c_{10}, 6)}$ are complex numbers such that $\lambda_{\dot{M}^1, \dot{M}^2}^{\dot M^3}=0$  if $N_{\dot{M}^1, \dot{M}^2}^{\dot{M}^3}=0$. We  now define an operator $Y$ on $\cu=L_{sl_2}(10, 0)\otimes L_{sl_2}(1, 0)\oplus L_{sl_2}(10, 6)\otimes L_{sl_2}(1, 0)$ as follows:
$$Y|_{{M}^1\times {M}^2}=\lambda_{\dot{M}^1, \dot{M}^2}^{\dot{L}(c_{10}, 0)}\Y_{{M}^1, {M}^2}^{L_{sl_2}(10, 0)}+\lambda_{\dot{M}^1, \dot{M}^2}^{\dot{L}(c_{10}, 6)}\Y_{{M}^1, {M}^2}^{L_{sl_2}(10, 6)}.$$
In the following we shall prove that $(\cu, Y)$ is a vertex operator algebra. First, we  prove the commutativity: For any $w^1\in M^1, w^2\in M^2, w^3\in M^3$ and $w'\in \cu '$,
\begin{align*}
E\langle w', Y(w^1, z_1)Y(w^2, z_2)w^3\rangle
=E\langle w', Y(w^2, z_2)Y(w^1, z_1)w^3\rangle.
\end{align*}
By the definition of $Y$, we have
\begin{align*}
E\langle& w', Y(w^1, z_1)Y(w^2, z_2)w^3\rangle\\
&=E\langle w',\lambda_{\dot{M}^1, \dot{L}_{sl_2}(10, 0)}^{\dot{L}_{sl_2}(10, 0)}\lambda_{\dot{M}^2, \dot{M}^3}^{\dot{L}_{sl_2}(10, 0)}\Y_{{M}^1, L_{sl_2}(10, 0)}^{L_{sl_2}(10, 0)}(w^1, z_1)\Y_{{M}^2, {M}^3}^{L_{sl_2}(10, 0)}(w^2, z_2)w^3\rangle\\
&\ \ +E\langle w',\lambda_{\dot{M}^1, \dot{L}_{sl_2}(10, 0)}^{\dot{L}_{sl_2}(10, 6)}\lambda_{\dot{M}^2, \dot{M}^3}^{\dot{L}_{sl_2}(10, 0)}\Y_{{M}^1, L_{sl_2}(10, 0)}^{L_{sl_2}(10, 6)}(w^1, z_1)\Y_{{M}^2, {M}^3}^{L_{sl_2}(10, 0)}(w^2, z_2)w^3\rangle\\
&\ \ +E\langle w',\lambda_{\dot{M}^1, \dot{L}_{sl_2}(10, 6)}^{\dot{L}_{sl_2}(10, 0)}\lambda_{\dot{M}^2, \dot{M}^3}^{\dot{L}_{sl_2}(10, 6)}\Y_{{M}^1, L_{sl_2}(10, 6)}^{L_{sl_2}(10, 0)}(w^1, z_1)\Y_{{M}^2, {M}^3}^{L_{sl_2}(10, 6)}(w^2, z_2)w^3\rangle\\
&\ \ +E\langle w', \lambda_{\dot{M}^1, \dot{L}_{sl_2}(10, 6)}^{\dot{L}_{sl_2}(10, 6)}\lambda_{\dot{M}^2, \dot{M}^3}^{\dot{L}_{sl_2}(10, 6)}\Y_{{M}^1, L_{sl_2}(10, 6)}^{L_{sl_2}(10, 6)}(w^1, z_1)\Y_{{M}^2, {M}^3}^{L_{sl_2}(10, 6)}(w^2, z_2)w^3\rangle,
\end{align*}
and
\begin{align*}
E\langle &w', Y(w^2, z_2)Y(w^1, z_1)w^3\rangle\\
&=E\langle w',\lambda_{\dot{M}^2, \dot{L}_{sl_2}(10, 0)}^{\dot{L}_{sl_2}(10, 0)}\lambda_{\dot{M}^1, \dot{M}^3}^{\dot{L}_{sl_2}(10, 0)}\Y_{{M}^2, L_{sl_2}(10, 0)}^{L_{sl_2}(10, 0)}(w^2, z_2)\Y_{{M}^1, {M}^3}^{L_{sl_2}(10, 0)}(w^1, z_1)w^3\rangle\\
&\ \ +E\langle w',\lambda_{\dot{M}^2, \dot{L}_{sl_2}(10, 0)}^{\dot{L}_{sl_2}(10, 6)}\lambda_{\dot{M}^1, \dot{M}^3}^{\dot{L}_{sl_2}(10, 0)}(w^2, z_2)\Y_{{M}^1, L_{sl_2}(10, 0)}^{L_{sl_2}(10, 6)}\Y_{{M}^1, {M}^3}^{L_{sl_2}(10, 0)}(w^1, z_2)w^3\rangle\\
&\ \ +E\langle w', \lambda_{\dot{M}^2, \dot{L}_{sl_2}(10, 6)}^{\dot{L}_{sl_2}(10, 0)}\lambda_{\dot{M}^1, \dot{M}^3}^{\dot{L}_{sl_2}(10, 6)}(w^2, z_2)\Y_{{M}^2, L_{sl_2}(10, 6)}^{L_{sl_2}(10, 0)}\Y_{{M}^1, {M}^3}^{L_{sl_2}(10, 6)}(w^1, z_2)w^3\rangle\\
&\ \ +E\langle w', \lambda_{\dot{M}^2, \dot{L}_{sl_2}(10, 6)}^{\dot{L}_{sl_2}(10, 6)}\lambda_{\dot{M}^1, \dot{M}^3}^{\dot{L}_{sl_2}(10, 6)}(w^2, z_2)\Y_{{M}^2, L_{sl_2}(10, 6)}^{L_{sl_2}(10, 6)}\Y_{{M}^1, {M}^3}^{L_{sl_2}(10, 6)}(w^1, z_2)w^3\rangle.
\end{align*}
Note that if $w^1\in \dot{M}^1, w^2\in \dot{M}^2, w^3\in \dot{M}^3$, then we have  \begin{align*}
E\langle w', Y(w^1, z_1)&Y(w^2, z_2)w^3\rangle
=E\langle w', Y(w^2, z_2)Y(w^1, z_1)w^3\rangle,
\end{align*}
implies
\begin{align*}
&E\langle w',\lambda_{\dot{M}^1, \dot{L}_{sl_2}(10, 0)}^{\dot{L}_{sl_2}(10, 0)}\lambda_{\dot{M}^2, \dot{M}^3}^{\dot{L}_{sl_2}(10, 0)}\Y_{\dot{M}^1, \dot{L}_{sl_2}(10, 0)}^{\dot{L}_{sl_2}(10, 0)}(w^1, z_1)\Y_{\dot{M}^2, \dot{M}^3}^{\dot{L}_{sl_2}(10, 0)}(w^2, z_2)w^3\rangle\\
&\ \ +E\langle w',\lambda_{\dot{M}^1, \dot{L}_{sl_2}(10, 0)}^{\dot{L}_{sl_2}(10, 6)}\lambda_{\dot{M}^2, {M}^3}^{\dot{L}_{sl_2}(10, 0)}\Y_{\dot{M}^1, \dot{L}_{sl_2}(10, 0)}^{\dot{L}_{sl_2}(10, 6)}(w^1, z_1)\Y_{\dot{M}^2, \dot{M}^3}^{\dot{L}_{sl_2}(10, 0)}(w^2, z_2)w^3\rangle\\
&\ \ +E\langle w',\lambda_{\dot{M}^1, \dot{L}_{sl_2}(10, 6)}^{\dot{L}_{sl_2}(10, 0)}\lambda_{\dot{M}^2, \dot{M}^3}^{\dot{L}_{sl_2}(10, 6)}\Y_{\dot{M}^1, \dot{L}_{sl_2}(10, 6)}^{\dot{L}_{sl_2}(10, 0)}(w^1, z_1)\Y_{\dot{M}^2, \dot{M}^3}^{\dot{L}_{sl_2}(10, 6)}(w^2, z_2)w^3\rangle\\
&\ \ +E\langle w',\lambda_{\dot{M}^1, \dot{L}_{sl_2}(10, 6)}^{\dot{L}_{sl_2}(10, 6)}\lambda_{\dot{M}^2, \dot{M}^3}^{\dot{L}_{sl_2}(10, 6)}\Y_{\dot{M}^1, \dot{L}_{sl_2}(10, 6)}^{\dot{L}_{sl_2}(10, 6)}(w^1, z_1)\Y_{\dot{M}^2, \dot{M}^3}^{\dot{L}_{sl_2}(10, 6)}(w^2, z_2)w^3\rangle\\
&=E\langle w',\lambda_{\dot{M}^2, \dot{L}_{sl_2}(10, 0)}^{\dot{L}_{sl_2}(10, 0)}\lambda_{\dot{M}^1, \dot{M}^3}^{\dot{L}_{sl_2}(10, 0)}\Y_{\dot{M}^2, \dot{L}_{sl_2}(10, 0)}^{\dot{L}_{sl_2}(10, 0)}(w^2, z_2)\Y_{\dot{M}^1, \dot{M}^3}^{\dot{L}_{sl_2}(10, 0)}(w^1, z_1)w^3\rangle\\
&\ \ +E\langle w',\lambda_{\dot{M}^2, \dot{L}_{sl_2}(10, 0)}^{\dot{L}_{sl_2}(10, 6)}\lambda_{\dot{M}^1, \dot{M}^3}^{\dot{L}_{sl_2}(10, 0)}(w^2, z_2)\Y_{{M}^1, \dot{L}_{sl_2}(10, 0)}^{\dot{L}_{sl_2}(10, 6)}\Y_{\dot{M}^1, \dot{M}^3}^{\dot{L}_{sl_2}(10, 0)}(w^1, z_2)w^3\rangle\\
&\ \ +E\langle w',\lambda_{\dot{M}^2, \dot{L}_{sl_2}(10, 6)}^{\dot{L}_{sl_2}(10, 0)}\lambda_{\dot{M}^1, \dot{M}^3}^{\dot{L}_{sl_2}(10, 6)}(w^2, z_2)\Y_{\dot{M}^2, \dot{L}_{sl_2}(10, 6)}^{\dot{L}_{sl_2}(10, 0)}\Y_{\dot{M}^1, \dot{M}^3}^{\dot{L}_{sl_2}(10, 6)}(w^1, z_2)w^3\rangle\\
&\ \ +E\langle w',\lambda_{\dot{M}^2, \dot{L}_{sl_2}(10, 6)}^{\dot{L}_{sl_2}(10, 6)}\lambda_{\dot{M}^1, \dot{M}^3}^{\dot{L}_{sl_2}(10, 6)}(w^2, z_2)\Y_{\dot{M}^2, \dot{L}_{sl_2}(10, 6)}^{\dot{L}_{sl_2}(10, 6)}\Y_{\dot{M}^1, \dot{M}^3}^{\dot{L}_{sl_2}(10, 6)}(w^1, z_2)w^3\rangle.
\end{align*}
Since for  $i\in\{0, 6\}$,
$$\{E\langle w',\Y_{\dot{M}^1, \dot{M}}^{\dot{L}_{sl_2}(10, i)}(w^1, z_1)\Y_{\dot{M}^2, \dot{M}^3}^{\dot{M}}(w^2, z_2)w^3\rangle|\dot{M}\}$$
 are linearly independent and
\begin{align*}
sp&an\{E\langle w',\Y_{{M}^1, M}^{{L}_{sl_2}(10, i)}(w^1, z_1)\Y_{{M}^2, {M}^3}^{M}(w^2, z_2)w^3\rangle|M\}\\
&=span\{E\langle w',\Y_{{M}^1, N}^{{L}_{sl_2}(10, i)}(w^1, z_1)\Y_{{M}^2, {M}^3}^{N}(w^2, z_2)w^3\rangle|N\},
\end{align*}
 we have \begin{align*}
E\langle w', Y(w^1, z_1)Y(w^2, z_2)w^3\rangle
=E\langle w', Y(w^2, z_2)Y(w^1, z_1)w^3\rangle
\end{align*}
holds for any $w^1\in M^1, w^2\in M^2, w^3\in M^3$ and $w'\in \cu '$. Similarly, we have
\begin{align*}
E\langle w', Y(w^1, z_1)Y(w^2, z_2)w^3\rangle
=E\langle w', Y(Y(w^1, z_1-z_2)w^2, z_2)w^3\rangle.
\end{align*}
holds for any $w^1\in M^1, w^2\in M^2, w^3\in M^3$ and $w'\in \cu '$. Thus, from Proposition 1.7 in \cite{H4}, we have $(\cu, Y)$ is a vertex operator algebra. Note that  $\cu$ is  an extension vertex operator algebra of $L_{sl_2}(10,0)\otimes L_{sl_2}(1,0)$, it follows from Theorem \ref{u1} that there is a unique vertex operator algebra structure on $\cu$, thus $\overline U$ must be isomorphic to the vertex operator algebra in Theorem \ref{vvoa5}. This finishes the proof.
\qed

\vskip0.5cm
 We now construct the vertex operator algebra realizing the modular invariant of type $( E_8, A_{30})$. Consider the conformal inclusions
 $$SU(2)_{28}\subset (G_2)_1\ \ \ \ \ \ \ \ (L_{sl_2}(28, 0)\subset L_{G_2}(1,0)).$$
  From the discussion in Section 3.3, we know that the decomposition of $ L_{G_2}(1,0)$ as $L_{sl_2}(28,0)$-module is as follows:
  \begin{align*}
  L_{G_2}(1,0)=L_{sl_2}(28, 0)\oplus L_{sl_2}(28, 10)\oplus L_{sl_2}(28, 18)\oplus L_{sl_2}(28, 28).
  \end{align*}
Consider the vertex operator algebra $L_{G_2}(1, 0) \otimes L_{sl_2}(1, 0)$, this vertex operator algebra contains $L_{sl_2}(28, 0)\otimes L_{sl_2}(1,0)$. The decomposition of $L_{G_2}(1, 0) \otimes L_{sl_2}(1, 0)$ as  $L_{sl_2}(28, 0)\otimes L_{sl_2}(1,0)$-module is as follows:
\begin{align*}
 L_{G_2}(1,0) \otimes L_{sl_2}(1, 0)&=L_{sl_2}(28, 0)\otimes L_{sl_2}(1, 0)\oplus L_{sl_2}(28, 10)\otimes L_{sl_2}(1, 0)\\
 &\oplus L_{sl_2}(28, 18)\otimes L_{sl_2}(1, 0)\oplus L_{sl_2}(28, 28)\otimes L_{sl_2}(1, 0).
 \end{align*}
By the discussion at the beginning of this subsection, we know that $L_{sl_2}(29, 0)$ is a vertex operator subalgebra of $L_{G_2}(1, 0) \otimes L_{sl_2}(1, 0)$, and we have
\begin{align*}
L_{G_2}(1,0) \otimes L_{sl_2}(1, 0)&=\mathop\oplus\limits_{\substack{
0\leq s\leq 29,\\
s\equiv 0 \mod  2}}
L(c_{28}, h^{28}_{1,s+1})\otimes  L_{sl_{2}}(29, s)\\
&\bigoplus\mathop\oplus\limits_{\substack{
0\leq s\leq 29,\\
s\equiv 0 \mod  2}}
L(c_{28}, h^{28}_{11,s+1})\otimes  L_{sl_{2}}(29, s)\\
&\bigoplus\mathop\oplus\limits_{\substack{
0\leq s\leq 29,\\
s\equiv 0 \mod  2}}
L(c_{28}, h^{28}_{19,s+1})\otimes  L_{sl_{2}}(29, s)\\
&\bigoplus\mathop\oplus\limits_{\substack{
0\leq s\leq 29,\\
s\equiv 0 \mod  2}}
L(c_{28}, h^{28}_{29,s+1})\otimes  L_{sl_{2}}(29, s).
\end{align*}
Thus the commutant vertex operator algebra of $L_{sl_2}(29, 0)$ in $L_{G_2}(1, 0) \otimes L_{sl_2}(1, 0)$ is $$U=L(c_{28}, 0)\oplus L(c_{28},h_{11,1}^{28})\oplus L(c_{28},h_{19,1}^{28})\oplus L(c_{28},h_{29,1}^{28}).$$
\begin{theorem}\label{11.1}
Let $U=L(c_{28}, 0)\oplus L(c_{28},h_{11,1}^{28})\oplus L(c_{28},h_{19,1}^{28})\oplus L(c_{28},h_{29,1}^{28})$ be the vertex operator algebra defined above. Then there is a unitary vertex operator algebra structure on $U$. Moreover, $U$  satisfies conditions (i) and (ii).
\end{theorem}
\pf This follows from an argument similar to that of Theorem \ref{vvoa5}.\qed
\vskip0.25cm
Therefore, by Theorem \ref{invariant}, $U=L(c_{28}, 0)\oplus L(c_{28},h_{11,1}^{28})\oplus L(c_{28},h_{19,1}^{28})\oplus L(c_{28},h_{29,1}^{28})$ realizes a modular invariant of $L(c_{28}, 0)$. Obviously, this modular invariant is equal to the  modular invariant of type $( E_8, A_{30})$:
\begin{tiny}
\begin{align*}
\frac{1}{2}\sum\limits_{p=1}^{30}&\bigg\{(Z_{1,q}(u,\tau_1)+Z_{11,q}(u,\tau_1)+Z_{19,q}(u,\tau_1)+Z_{29,q}(u,\tau_1))\overline{(Z_{1,q}(v,\tau_2)+Z_{11,q}(v,\tau_2)+Z_{19,q}(v,\tau_2)+Z_{29,q}(v,\tau_2))}\\ &(Z_{7,q}(u,\tau_1)+Z_{13,q}(u,\tau_1)+Z_{17,q}(u,\tau_1)+Z_{23,q}(u,\tau_1))\overline{(Z_{7,q}(v,\tau_2)+Z_{13,q}(v,\tau_2)+Z_{17,q}(v,\tau_2)+Z_{23,q}(v,\tau_2))}\bigg\}.
\end{align*}
\end{tiny}
From the modular invariant we immediately know that the vector spaces
\begin{align*}
&L(c_{28}, h_{1,q}^{28})\oplus L(c_{28},h_{11,q}^{28})\oplus L(c_{28},h_{19,q}^{28})\oplus L(c_{28},h_{29,q}^{28}),\\
&L(c_{28}, h_{7,q}^{28})\oplus L(c_{28},h_{13,q}^{28})\oplus L(c_{28},h_{17,q}^{28})\oplus L(c_{28},h_{23,q}^{28}),\ 1\leq q \leq 16,
  \end{align*} have irreducible $U$-module structure and these are exactly all the irreducible $U$-modules.

A similar argument as  in the proof of Theorem \ref{vvoa4} gives the uniqueness of the vertex operator algebra structure on  $L(c_{28}, 0)\oplus L(c_{28},h_{11,1}^{28})\oplus L(c_{28},h_{19,1}^{28})\oplus L(c_{28},h_{29,1}^{28})$.
\begin{theorem}
Let $\overline U$ be an extension vertex operator algebra of $L(c_{28}, 0)$ such that $\overline U$ is isomorphic to $L(c_{28}, 0)\oplus L(c_{28},h_{11,1}^{28})\oplus L(c_{28},h_{19,1}^{28})\oplus L(c_{28},h_{29,1}^{28})$ as $L(c_{28}, 0)$-module. Then $\overline U$ is isomorphic to the vertex operator algebra defined above.
\end{theorem}
\subsection{Realization of modular invariants of type $(A_{10}, E_6)$ and $(A_{28}, E_8)$}
In this subsection we  construct extension vertex operator algebras of $L(c_m,0)$ which realize the modular invariants of type $(A_{10}, E_6)$ and $(A_{28}, E_8)$.
Consider the vertex operator algebra $L_{sl_{2}}(9,0)\otimes  L_{sl_{2}}(1, 0)$. By the discussion in subsection 4.2, the decomposition of  $L_{sl_{2}}(9,0)\otimes  L_{sl_{2}}(1, 0)$  as $L_{sl_{2}}(10, 0)$-module is as follows:
\begin{align*}
L_{sl_2}(9,0) \otimes L_{sl_2}(1, 0)&=\mathop\oplus\limits_{\substack{
0\leq s\leq 10,\\
s\equiv 0 \mod  2}}
L(c_{9}, h^{9}_{1,s+1})\otimes  L_{sl_{2}}(10, s).
\end{align*}
Since there is a vertex operator algebra structure on $L_{sl_{2}}(10, 0)\oplus L_{sl_{2}}(10, 6)$, it follows from Theorem \ref{kvoa3} that there is a vertex operator algebra structure on $L(c_9, 0)\oplus L(c_9, h_{1,7}^9)$.
\begin{theorem}\label{kvoa15}
Let $U=L(c_9, 0)\oplus L(c_9, h_{1,7}^9)$ be  the vertex operator algebra defined above. Then $U$  satisfies conditions (i) and (ii).
\end{theorem}
\pf Note that $L_{B_2}(1,0)$ is a simple vertex operator algebra, then $U$ is a simple vertex operator algebra by Theorem \ref{kvoa3}. Since $L(c_{9}, 0)$ satisfies conditions (i) and (ii), it follows from Theorems \ref{C 2}, \ref{cvoa2} that $U=L(c_9, 0)\oplus L(c_9, h_{1,7}^9)$ is rational and $C_2$-cofinite. By a similar argument as in the proof of  Theorem \ref{vvoa5},  we could prove that $U$ is self dual. The proof is complete. \qed
\vskip0.25cm
Therefore, by Theorem \ref{invariant},  $U=L(c_9, 0)\oplus L(c_9, h_{1,7}^9)$ realizes a modular invariant of $L(c_9, 0)$. Obviously, this modular invariant is equal to the modular invariant of type $(A_{10}, E_6)$:
\begin{tiny}
\begin{align*}
\frac{1}{2}\sum\limits_{p=1}^{10}&\bigg\{(Z_{p,1}(u,\tau_1)+Z_{p,7}(u,\tau_1))\overline{(Z_{p,1}(v,\tau_2)+Z_{p,7}(v,\tau_2))}+(Z_{p,4}(u,\tau_1)+Z_{p,8}(u,\tau_1))\overline{(Z_{p,4}(v,\tau_2)+Z_{p,8}(v,\tau_2))}\\ &+(Z_{p,5}(u,\tau_1)+Z_{p,11}(u,\tau_1))\overline{(Z_{p,5}(v,\tau_2)+Z_{p,11}(v,\tau_2))}\bigg\}.
\end{align*}
\end{tiny}
From the modular invariant we immediately know that the vector spaces
\begin{align*}
&L(c_9, h_{p,1}^9)\oplus L(c_9, h_{p,7}^9),\\
 &L(c_9, h_{p,4}^9)\oplus L(c_9, h_{p,8}^9),\\
 &L(c_9, h_{p,5}^9)\oplus L(c_9, h_{p,11}^9),\ 1\leq p\leq 6,
 \end{align*}  have irreducible $U$-module structures and these are exactly all the irreducible $U$-modules.
Moreover, we have
\begin{theorem}\label{AE_6}
Let $\overline U$ be an extension vertex operator algebra of $L(c_{9}, 0)$ such that $\overline U$ is isomorphic to $L(c_9, 0)\oplus L(c_9, h_{1,7}^9)$ as $L(c_{9}, 0)$-module. Then $\overline U$ is isomorphic to the vertex operator algebra defined above.
\end{theorem}
\pf Assume that $\overline U=L(c_{9}, 0)\oplus L(c_{9},h_{7,1}^{10})$ is an extension vertex operator algebra of $L(c_{9}, 0)$, it follows from Theorems \ref{kvoa2}, \ref{kvoa3} that there is a vertex operator algebra structure on $L_{sl_2}(10, 0)\oplus L_{sl_2}(10, 6)$ such that $L_{sl_2}(10, 0)\oplus L_{sl_2}(10, 6)$ is an extension vertex operator algebra  of $L_{sl_2}(10, 0)$. Since there is a unique vertex operator algebra structure on $L_{sl_2}(10, 0)\oplus L_{sl_2}(10, 6)$,  $\overline U$ must be  isomorphic to the vertex operator algebra defined above. The proof is complete.
\qed

\vskip0.5cm
We now construct the vertex operator algebra realizing the modular invariant of type $(A_{28}, E_8)$. Consider the vertex operator algebra  $L_{sl_{2}}(27,0)\otimes  L_{sl_{2}}(1, 0)$. By the discussion in subsection 4.2, the decomposition of  $L_{sl_{2}}(27,0)\otimes  L_{sl_{2}}(1, 0)$  as $L_{sl_{2}}(28, 0)$-module is as follows:
\begin{align*}
L_{sl_2}(27,0) \otimes L_{sl_2}(1, 0)&=\mathop\oplus\limits_{\substack{
0\leq s\leq 28,\\
s\equiv 0 \mod  2}}
L(c_{27}, h^{27}_{1,s+1})\otimes  L_{sl_{2}}(28, s).
\end{align*}
Since there is a vertex operator algebra structure on $$L_{sl_{2}}(28, 0)\oplus L_{sl_{2}}(28, 10)\oplus L_{sl_{2}}(28, 18)\oplus L_{sl_{2}}(28, 28),$$ it follows from Theorem \ref{kvoa3} that there is a vertex operator algebra structure on $$L(c_{27}, 0)\oplus L(c_{27},h_{1,11}^{27})\oplus L(c_{27},h_{1,19}^{27})\oplus L(c_{27},h_{1,29}^{27})$$ such that $L(c_{27}, 0)$ is a vertex operator subalgebra of this vertex operator algebra.

\begin{theorem}\label{kvoa16}
Let $U=L(c_{27}, 0)\oplus L(c_{27},h_{1,11}^{27})\oplus L(c_{27},h_{1,19}^{27})\oplus L(c_{27},h_{1,29}^{27})$ be the vertex operator algebra defined above. Then $U$ satisfies conditions (1) and (2).
\end{theorem}
\pf This follows from an argument similar to that of Theorem \ref{kvoa15}.\qed
\vskip0.25cm
 Therefore, by Theorem \ref{invariant},  $U=L(c_{27}, 0)\oplus L(c_{27},h_{1,11}^{27})\oplus L(c_{27},h_{1,19}^{27})\oplus L(c_{27},h_{1,29}^{27})$ realizes a modular invariant  of  $L(c_{27}, 0)$. Obviously, this modular invariant is equal to the  modular invariant of type $(A_{28}, E_8)$:
\begin{tiny}
\begin{align*}
\frac{1}{2}\sum\limits_{p=1}^{28}& \bigg\{(Z_{p,1}(u,\tau_1)+Z_{p,11}(u,\tau_1)+Z_{p,19}(u,\tau_1)+Z_{p,29}(u,\tau_1))\overline{(Z_{p,1}(v,\tau_2)+Z_{p,11}(v,\tau_2)+Z_{p,19}(v,\tau_2)+Z_{p,29}(v,\tau_2))}\\
&(Z_{p,7}(u,\tau_1)+Z_{p,13}(u,\tau_1)+Z_{p,17}(u,\tau_1)+Z_{p,23}(u,\tau_1))\overline{(Z_{p,7}(v,\tau_2)+Z_{p,13}(v,\tau_2)+Z_{p,17}(v,\tau_2)+Z_{p,23}(v,\tau_2))} \bigg\}.
\end{align*}
\end{tiny}
From the modular invariant we immediately know that the vector spaces
\begin{align*}
&L(c_{27}, h_{p,1}^{27})\oplus L(c_{27},h_{p,11}^{27})\oplus L(c_{27},h_{p,19}^{27})\oplus L(c_{27},h_{p,29}^{27}),\\
&L(c_{27}, h_{p,7}^{27})\oplus L(c_{27},h_{p,13}^{27})\oplus L(c_{27},h_{p,17}^{27})\oplus L(c_{27},h_{p,23}^{27}),\ 1\leq p \leq 15,
  \end{align*} have irreducible $U$-module structures and these are exactly all the irreducible $U$-modules.
Moreover, 
using a similar argument as in the proof of  Theorem \ref{AE_6}, we have
\begin{theorem}\label{vvoa6}
Let $\overline U$ be an extension vertex operator algebra of $L(c_{27}, 0)$ such that $\overline U$ is isomorphic to $L(c_{27}, 0)\oplus L(c_{27},h_{1,11}^{27})\oplus L(c_{27},h_{1,19}^{27})\oplus L(c_{27},h_{1,29}^{27})$ as $L(c_{27}, 0)$-module. Then $\overline U$ is isomorphic to the vertex operator algebra defined above.
\end{theorem}
\subsection{Classification of preunitary vertex operator algebras with central charge $c<1$}
In this subsection we shall classify preunitary vertex operator algebras with central charge $c<1$.

Recall that a vertex operator algebra $V$ is preunitary if it is a unitary module for the Virasoro algebra. By a similar argument as in the proof of  Lemma 2.5 in  \cite{DLin}, we can prove the following:
\begin{proposition}\label{extension}
Let $V$ be a preunitary vertex operator algebra with central charge
$c<1$. Then the vertex operator subalgebra $\langle \w\rangle$ which is
generated by the conformal vector $\w$ is isomorphic to the vertex operator algebra $L(c_m, 0)$ for some integer $m\geq 2$. In particular, $V$ is an extension vertex operator algebra of $L(c_m, 0)$.
\end{proposition}
Now we can prove the second main result in this paper.
\begin{theorem}\label{kvoa14}
Let $V$ be a simple preunitary vertex operator algebra with central charge $c<1$. Assume that $V$ is of $CFT$ type  and  self dual, then $V$ is isomorphic to either $L(c_m,0)$ for some $m\geq 0$ or one of vertex operator algebras in Theorems \ref{vvoa1}, \ref{vvoa2}, \ref{vvoa5}, \ref{11.1}, \ref{kvoa15}, \ref{kvoa16}.
\end{theorem}
\pf From Proposition \ref{extension}, $V$ is an extension vertex operator algebra of $L(c_m, 0)$. Then $V$ satisfies conditions (i) and (ii). Thus  there is a modular invariant of $L(c_m, 0)$ associated to $V$ by Theorem \ref{invariant}. We first prove that this modular invariant is not equal to the modular invariants of type $(D_{2n+1}, A_{4n})$. Otherwise, if the modular invariant is equal to the modular invariants of type $(D_{2n+1}, A_{4n})$, then $V$ must be isomorphic to $L(c_m, 0)$, but   the modular invariants of type $(D_{2n+1}, A_{4n})$ is not equal to the modular invariant of type $(A_{m-1}, A_m)$, which is the modular invariant associated to the extension vertex operator algebra $L(c_m, 0)$, this is a contradiction. Similarly, the modular invariant is not equal to the modular invariants of type $(A_{4n-2}, D_{2n+1})$. We now prove that the modular invariant is not equal to the modular invariants of type $(E_{7}, A_{18})$. Otherwise, if the modular invariant is equal to the modular invariant of type $(E_{7}, A_{18})$, then $V$ is isomorphic to $L(c_{16}, 0)\oplus L(c_{16}, h^{16}_{17, 1})$ as $L(c_{16}, 0)$-module, note that $L(c_{16}, h^{16}_{17, 1})$ is a simple current $L(c_{16}, 0)$-module, it follows from Theorem \ref{simple current} that the vertex operator algebra structure on $L_{sl_2}(16, 0)\oplus L_{sl_2}(16, 16)$ is unique and is isomorphic to the vertex operator algebra in Theorem \ref{vvoa1}, but the modular invariant of type $(E_{7}, A_{18})$ is not equal to the modular invariant of type $(D_{10}, A_{18})$, this is a contradiction. Similarly, we could prove that the modular invariant is not equal to the modular invariant of type $(A_{16}, E_7)$. So the modular invariant is equal to one of the other modular invariants in Theorem \ref{ciz2}. By the discussion above we have  the vertex operator algebra realizing the modular invariant is unique.  Thus $V$ is isomorphic to either $L(c_m,0)$ or  one of vertex operator algebras in Theorems \ref{vvoa1}, \ref{vvoa2}, \ref{vvoa5}, \ref{11.1}, \ref{kvoa15}, \ref{kvoa16}. Now we prove that these vertex operator algebras are preunitary. Note that every vertex operator algebra here is an extension of
$L(c_m,0)$ and each irreducible $L(c_m,0)$-module is a unitary module for the Virasoro algebra. So all these
vertex operator algebras are preunitary. The proof is complete.
\qed
\begin{remark}
We have already proved the preunitary vertex operator algebras given in Theorem \ref{kvoa14} are unitary except
$L(c_9, 0)\oplus L(c_9, h_{1,7}^9)$ and $L(c_{27}, 0)\oplus L(c_{27},h_{1,11}^{27})\oplus L(c_{27},h_{1,19}^{27})\oplus L(c_{27},h_{1,29}^{27}).$ We believe that these two vertex operator algebras are also unitary vertex operator algebras, but we can not prove it in this paper.
\end{remark}

\end{document}